\documentclass[preprint, 12pt, number]{elsarticle}

\usepackage{changepage} 
\usepackage[mathscr]{euscript}
\usepackage{amssymb}
\usepackage{amsmath}
\usepackage{stackengine}
\usepackage{hyperref}
\usepackage{esint}
\usepackage{graphicx}
\usepackage{titlesec}
\usepackage{multirow}
\usepackage{multirow}
\usepackage{latexsym}
\usepackage{subcaption}
\usepackage{caption}
\usepackage{lipsum} 
\stackMath
\setcitestyle{square}
\usepackage{accents}
\usepackage{xcolor}
\usepackage{epstopdf}
\usepackage{bm}
\usepackage{float}
\usepackage[a4paper, margin=1in]{geometry} 

\begin{document}

\begin{frontmatter}

\title{The isogeometric boundary element algorithm for solving the plane strain problem of an elastic matrix containing an open material surface of arbitrary shape}

\author[umn]{Rohit Satish Patil}
\author[dhu]{Zhilin Han}
\author[umn]{Sofia G. Mogilevskaya\corref{cor1}}
\ead{mogil003@umn.edu}
\cortext[cor1]{Corresponding author.}

\address[umn]{Department of Civil, Environmental, and Geo-Engineering, University of Minnesota, 500 Pillsbury Drive S.E, 55414, Minneapolis, MN, USA}
\address[dhu]{College of Physics, Donghua University, Shanghai, 201620, China}

\begin{abstract}
 
 The paper presents the Isogeometric Boundary Element Method (IGABEM) algorithm for solving the plane strain problem of an isotropic linearly elastic matrix containing an open material surface of arbitrary shape. Theoretical developments are based on the use of the Gurtin-Murdoch model of material surfaces. The governing equations and the boundary conditions for the problem are reviewed, and analytical integral representations for the elastic fields everywhere in the material system are presented in terms of unknown traction jumps across the surface. To find the jumps, the problem is reduced to a system of singular boundary integral equations in terms of two unknown scalar components of the surface stress tensor. The system is solved numerically using the developed IGABEM algorithm in which NURBS are used to approximate the unknowns. The main steps of the algorithm are discussed and convergence studies are performed. The algorithm is validated using two benchmark problems involving the matrix subjected to a uniform far-field load and containing a surface along (i) a straight segment and (ii) a circular arc. Numerical examples are presented to illustrate the influence of governing parameters with a focus on the influence of curvature variation.


\end{abstract}

\end{frontmatter}

\section{Introduction}\label{Section:1}

Mechanical behavior of materials containing ultra-thin inhomogeneities/layers is of significant interest in engineering and materials science, particularly when the latter are stiff and subject to prestress. Although the study of such materials has a long history in continuum mechanics, most of the early works treated inhomogeneities/layers as rigid line inclusions, first considering single straight line inclusion and later extending the models to include multiple reinforcements. Comprehensive reviews of theoretical, numerical, and experimental investigations of this nature are available in \cite{Pingle2008, Jobin2019, Goudarzi2021, Patil2024}.

Almost at the same time, different models were developed that treated thin inhomogeneities/layers as elastic lines across which the displacements and tractions underwent jumps. The jump conditions, obtained by various asymptotic methods (see, i.e., the reviews of early work in \cite{Klarbring1991, Klarbring1998, Mishuris2004}), contained information on the thicknesses and elastic properties of the reinforcements. According to the classification of \cite{Benveniste2001}, stiff inhomogeneities/layers can be treated as membrane- or shell-type interfaces. Extensive reviews of the relevant literature can be found in, i.e., \cite{Benveniste2001,Rubin2004,Rizzoni2013,Baranova2022,Kushch2024} and the references therein. However, the above models generally did not include the effects due to prestress, restricting their applicability to more complex or realistic scenarios.

More recently, it has been suggested to model thin and stiff elastic inhomogeneities/layers as material surfaces using the Gurtin–Murdoch \cite{Gurtin1975,Gurtin1978} or Steigmann-Ogden theories \cite{Steigmann1997,Steigmann1999}, see the reviews in, i.e., \cite{Chhapadia2011, Firooz2021,Wang2011,Eremeyev2016_AM,Javili2018,Mogilevskaya2021_AMR}. In these theories, inhomogeneities/layers were treated as elastic membranes or shells of vanishing thickness, characterized by their own distinct mechanical properties and residual surface tension. The jump conditions across the surfaces included the conditions of continuity of displacements and jumps in tractions. These models extended classical elasticity by incorporating surface energy and elasticity effects, thereby providing a more complete description of the thin, stiff, and prestressed inhomogeneities/layers.

The developments of analytical and semi-analytical solutions for problems with material surfaces were initially restricted to materials reinforced by thinly coated fibers or particles, with coatings treated as Gurtin-Murdoch material surfaces. Analytical solutions based on the complete Gurtin–Murdoch model for such problems were reported in, i.e., \cite{He2006,Lim2006} for the surface along the boundary of a spherical cavity.
Semi-analytical solutions for materials with multiple circular fibers with Gurtin-Murdoch material surfaces were developed in \cite{Mogilevskaya2008,Jammes2009} and for materials with spherical particles in \cite{Kushch2011}.

The first solutions to the problems involving the Steigmann–Ogden materials surfaces along the boundaries of a single spherical particle or a circular fiber were obtained in \cite{Zemlyanova2018_JAM,Zemlyanova2018_IJSS}. Those solutions were further generalized for the case of multiple circular fibers in \cite{Han2018}. In \cite{Mogilevskaya2019}, analytical displacement representations were proposed for problems involving spherical and circular material surfaces described by the complete Gurtin-Murdoch and Steigmann–Ogden models.

Problems involving closed material surfaces of arbitrary shapes were also solved numerically. The Finite Element Method (FEM) solutions for these problems were reported in \cite{McBride2012,Javili2014,Chen2020,He2020}, whereas the Boundary Element Method (BEM) solutions in \cite{Dong2011,Dong2012,Xu2016}. The extensive reviews of the relevant numerical results can be found in \cite{Han2019,Chen2020,Firooz2021,Han2021_IJSS} and the references therein. 

As mentioned above, the Gurtin-Murdoch and Steigmann-Ogden models were mostly used for materials with coated fibers and particles, and therefore only closed material surfaces were considered. Starting in the 2020s, it was proposed to use the Gurtin-Murdoch and Steigmann-Ogden theories to model composite materials
reinforced by two-dimensional flexible membranes and shells. Such two-dimensional reinforcements can simulate ultrathin graphene nanoplatelets or graphene-oxide sheets used in the emerging generation of composites. The models of such kind require the concept of open material surfaces, i.e., surfaces that possess boundary curves or tip points.

The solution to the problem of an elastic matrix containing a single Gurtin–Murdoch material surface along a straight segment was first derived in anti-plane setting in \cite{Baranova2020_JE}; the solution to the corresponding plane strain problem was reported in \cite{Mogilevskaya2021_CST}. This work was extended in \cite{Patil2024} to include multiple straight material surfaces, in order to capture interaction effects. 
In \cite{Zemlyanova2023_PhysicaD}, the solution for the plane strain problem of an elastic matrix with a single straight material surface described by the Steigmann–Ogden theory was developed. The plane strain problem of an elastic bimaterial plane containing a Steigmann–Ogden material surface along a finite segment of the bimaterial interface was solved in \cite{Zemlyanova2025}. The solutions to problems involving Gurtin–Murdoch or Steigmann–Ogden material surfaces located along a single or multiple circular arcs were derived in \cite{Han2024_IJSS,Han2024_IJES,Han2024_JAM}.
In all of those studies, the analysis was restricted to straight-line and circular-arc geometries. That allowed for the expansions of the components of the surface stress tensor in either a series of Chebyshev polynomials of the second kind or in a series of trigonometric functions. To take care of the singularities near the tips, the series expansions were multiplied by the square-root weight functions. Such representations allowed an accurate evaluation of the elastic fields and stress-intensity factors. 

However, series expansions do not work for problems involving open material surfaces of irregular shapes; these problems can only be solved numerically. To date, only one FEM-based algorithm has been reported for open Gurtin–Murdoch material surfaces, and that was done in the antiplane setting \cite{Herrera-Garrido2025}. In that paper, the FEM formulation was developed for problems involving both open and closed material surfaces that could possess corner or tip points. The formulation was subsequently extended to account for the interactions between two surfaces in \cite{Herrera-GarridoMAandMogilevskaya2025}.

In this paper, we develop a novel two-dimensional BEM-based algorithm to solve the plane strain problem of an infinite matrix containing an open Gurtin–Murdoch material surface of arbitrary but sufficiently smooth shape. We adopt the concept of isogeometric analysis that was first used by Hughes et al. \cite{Hughes2005, Cottrell2009} in the context of FEM. In the analysis, the same Non-Uniform Rational B-Splines functions (NURBS) are used to represent curved geometries and approximate field variables. The concept allows for exact representation of geometries and for higher order inter-element continuity for the functions involved, which gives the isogeometric FEM computational advantages over standard, polynomial-based FEM. 
The isogeometric FEM is now a well-established method that has been used in numerous engineering applications, see, e.g., \cite{Pivovarov2025, Su2024, Gupta2023}, for extensive literature reviews.

The advantages of isogeometric analysis and BEM were combined in the Isogeometric Boundary Element Method in \cite{Simpson2012, Simpson2013}. IGABEM has been successfully applied to solve steady state potential problems, as well as those of elasticity, viscous flow, acoustics, etc. A review of the latest developments in IGABEM and its various engineering applications can be found in \cite{Beer2020,Beer2025}. Those also include applications in fracture mechanics, e.g., \cite{Rocha2024}, \cite{Andrade2023}.

Here, we propose for the first time using IGABEM for modeling problems with open material surfaces of varying curvature. Our choice is motivated by the need to accurately represent the geometry of the material surface and to guaranty that the approximations of the surface stress tensor components, involved in the governing integral representations, possess required high order smoothness conditions. IGABEM is ideally suited to address these needs and, therefore is used here as the modeling tool.
 
The structure of the paper is as follows. The problem formulation is given in Section~\ref{Section:2}, while the governing equations of the Gurtin–Murdoch surface elasticity model in the plane strain setting are reviewed in Section~\ref{Section:3}. The integral representations of the elastic fields are presented in Section~\ref{Section:4} with the details provided in \ref{Section:Appendix A}. In Section~\ref{Section:5}, the basic steps of the proposed IGABEM algorithm are presented, including reviews of the theories of B-splines and NURBS, approximations of the geometry and surface stress components, the boundary integral equations in parametric forms, and their numerical solutions. Section~\ref{Sub-Section:6.1} presents several numerical examples designed to validate the proposed IGABEM algorithm using the available benchmark solutions and to illustrate the influence of curvature variations on elastic fields. The convergence study of the algorithm is reported in \ref{Section:Appendix B}. Finally, the concluding remarks are presented in Section~\ref{Section:7}.
 
\section{Problem Formulation}\label{Section:2}

\begin{figure}[htb]
	\centering
	\includegraphics[width=0.75\textwidth]{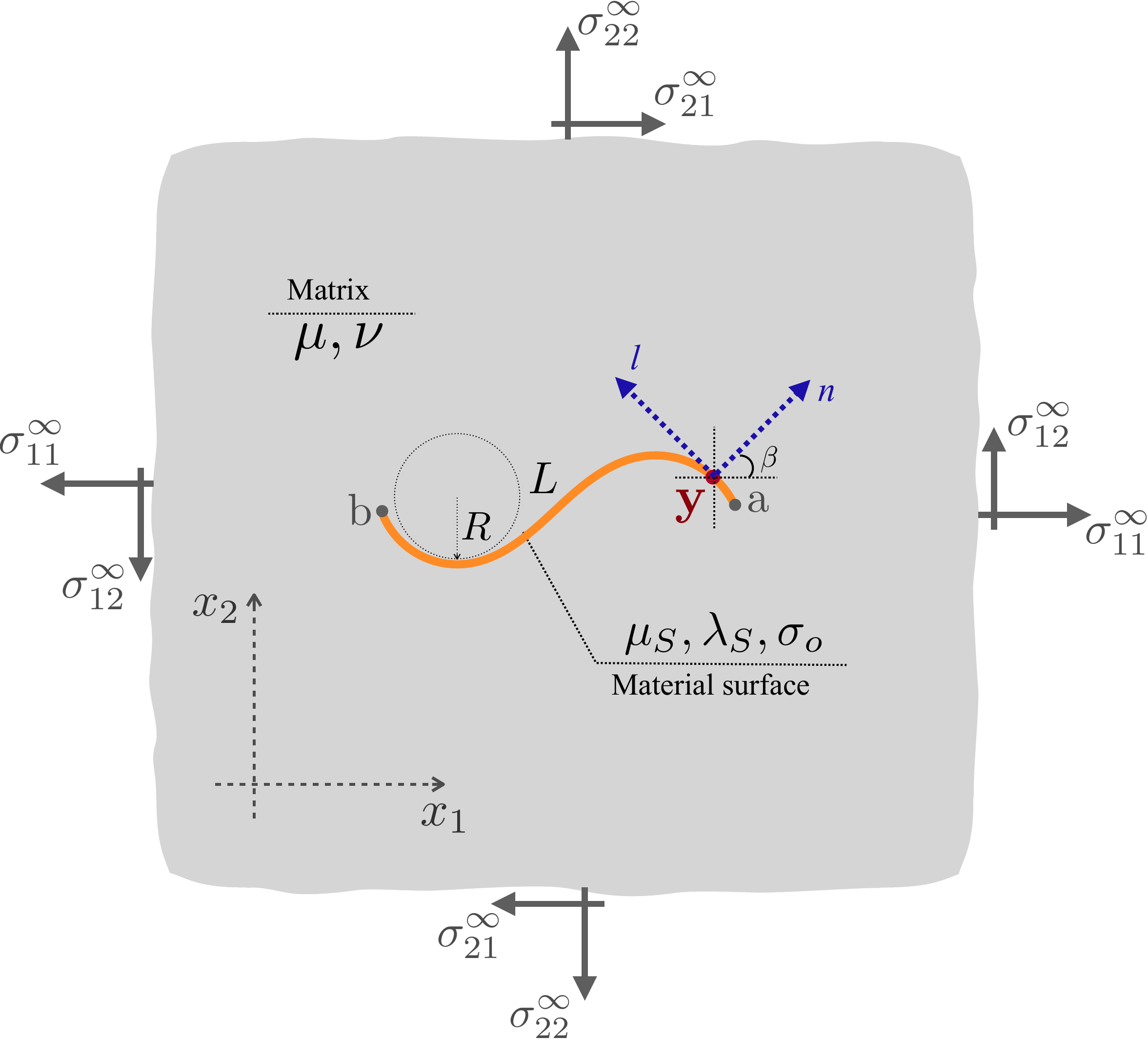}
	\caption{Problem formulation.}
    \label{Fig:1}
\end{figure}

Consider the two-dimensional plane strain problem involving an infinite isotropic linearly elastic matrix subjected to uniform far-field stress $\boldsymbol{\sigma}^{\infty}$ with components $\sigma^{\infty}_{11}$, $\sigma^{\infty}_{12}$, $\sigma^{\infty}_{22}$. The matrix contains an open material surface (curve, in two-dimensions) $L$ of an arbitrary but sufficiently smooth shape with tip points $\mathbf{a}$, $\mathbf{b}$, see Fig.~\ref{Fig:1}. The Gurtin-Murdoch theory of material surfaces is adopted, in which the surface is treated as a membrane of negligible thickness attached to the bulk without slipping. The mechanical properties of the membrane are defined by its surface shear modulus $\mu_{S}$, surface Lamé parameter $\lambda_{S}$, and surface tension $\sigma_{0}$, each having dimension of $\mathrm{N/m}$, while the isotropic elastic matrix is described by its shear modulus $\mu$, which have dimensions of $\mathrm{N/m^{2}}$, and Poisson's ratio $\nu$.
The goal is to design an efficient IGABEM algorithm for accurate evaluation of the elastic fields at any given point within the material system.
 
\section{Review of governing equations of the Gurtin-Murdoch model for the plane strain case}\label{Section:3}

The governing equations for the Gurtin-Murdoch model were derived in \cite{Gurtin1975}, \cite{Gurtin1978} and particularized for the plane strain case in \cite{Mogilevskaya2021_CST}, \cite{Mogilevskaya2021_AMR}. Here, we review the latter equations for the case of an open material surface $L$ of arbitrary sufficiently smooth shape. The model comprises the standard Navier equation that describes the bulk displacement field, complemented by the following conditions that describe the behavior of displacements and tractions across the surface:

\begin{equation}\label{Eq:1}
    \begin{gathered}
        u_{1}^{+}= u_{1}^{-}= u_{1}\:, \\[6pt]
        u_{2}^{+}= u_{2}^{-}= u_{2}\:,
    \end{gathered}
\end{equation}

\begin{equation}\label{Eq:2}
    \begin{gathered}
         \Delta t_{l} = t_{l}^{+} - t_{l}^{-} =  \frac{\partial \sigma^{S}}{\partial s} + \sigma_{0} \frac{\omega^{S}}{R}\:, \\[10pt]
        \Delta t_{n} = t_{n}^{+} - t_{n}^{-} = - \frac{\sigma^{S}}{R} + \sigma_{0} \frac{\partial\omega^{S}}{\partial s}\:.
    \end{gathered}
\end{equation}                      
The superscripts   ``+"  and ``--", used in Eqs.~(\ref{Eq:1})-(\ref{Eq:2}) and throughout the paper, represent the limit values of the corresponding fields as the surface is approached from the direction of the normal vector or from the opposite direction, respectively. The quantities $u_{1}$ and $u_{2}$ of Eq.~(\ref{Eq:1}) denote the components of the displacement vector in the Cartesian coordinates ($x_{1}, x_{2}$), while $t_{n}$ and $t_{l}$ of Eq.~(\ref{Eq:2}) represent the normal and shear components of the bulk tractions in the local coordinates ($n, l$), see Fig.~\ref{Fig:1}. Additional quantities involved in Eq.~(\ref{Eq:2}) are: the arc length parameter $s$ (measured from the tip $\bold{a}$), the local radius of curvature $R=R(s)$, and the only non-vanishing components $\sigma^{S}$ and $\omega^{S}$ of the surface stress tensor that are defined as 
\begin{equation}\label{Eq:3}
    \begin{gathered}
        \sigma^{S} = \sigma_{0} + (\lambda_{S} + 2\mu_{S}) \varepsilon^{S}\:, \\[10pt]
        \omega^{S} = - \frac{u_{l}}{R} + \frac{\partial u_{n}}{\partial s}\:,
    \end{gathered}
\end{equation}
in which $u_{n}$ and $u_{l}$ are the components of displacements in the local coordinates ($n, l$), and $\varepsilon^{S}$ is the only non-vanishing component of the surface strain tensor that can be expressed via the displacements as 
\begin{equation}\label{Eq:4}
 \varepsilon^{S} = \frac{u_{n}}{R} + \frac{\partial u_{l}}{\partial s}\:. \\[10pt]
\end{equation}
In addition, the following conditions must be enforced at the tips $\mathbf{a}$, $\mathbf{b}$ of the surface:
\begin{equation}\label{Eq:5}
    \begin{gathered}
         \sigma^{S}=0, \\[5 pt]
         \sigma_{0} \omega^{S} = 0\:. 
    \end{gathered}
\end{equation}

\section{Governing integral representations}\label{Section:4}
If the displacements in the matrix are expressed in the form of a single-layer elastic potential, \cite{Kupradze1965}, the conditions of continuity of the displacements and jumps in tractions are automatically fulfilled, as explained in \cite{Mogilevskaya2014}, \cite{CrouchMogilevskaya2024}. 
The use of a single layer elastic potential leads to the following expression for each Cartesian displacement component defined at the point $\mathbf{x}$ located everywhere in the domain, including at $L$ (as the single layer potential is continuous across $L$), via the traction jumps:
\begin{equation}\label{Eq:6}
    \begin{gathered}
u_{k}(\mathbf{x}) = u^{\infty}_{k} (\mathbf{x}) + \int_{L} \Delta t_{j} (\mathbf{y}) G_{kj} (\mathbf{x}, \mathbf{y})\, \mathrm{d}s_{\mathbf{y}}\:;\quad k,j = 1,2\:,
    \end{gathered}
\end{equation}
where $s_{\mathbf{y}}$ is the arc length of the material surface at the point $\mathbf{y}\in L$, $u^{\infty}_{k} (\mathbf{x})$ is the $k$-th component of the displacements at the point $\mathbf{x}$, located anywhere in the homogeneous plane (without material surface), due to the far-field load $\bm{\sigma}^{\infty}$, $\Delta t_{j} (\mathbf{y})= t_{j}^{+} (\mathbf{y}) - t_{j}^{-} (\mathbf{y})$ is the $j$-th Cartesian component of the traction jump across the surface at point $\mathbf{y}$, and the repeated index implies summation. The kernel $G_{kj} (\mathbf{x},\mathbf{y})$ of Eq.~(\ref{Eq:7}) is the following Kelvin fundamental solution:
\begin{equation}\label{Eq:7}
    \begin{gathered}
    G_{kj} (\mathbf{x}, \mathbf{y})= \frac{1}{2 \mathrm{\pi} \mu (\kappa +1)} [- \kappa \,\delta_{kj} \operatorname{ln} r + r,_{k} r,_{j}]\:,
    \end{gathered}
\end{equation}
in which $\delta_{kj}$ is Kronecker's symbol, $\kappa=3-4\nu$, $r = | \mathbf{x} -  \mathbf{y} |$, $r_{,k} = {\partial r}/{\partial x_{k}}$.

Using a standard coordinate transformation procedure and Eq.~(\ref{Eq:6}), the integral representations of the displacement components $u_{l},\, u_{n}$ in the local coordinate system can be found in terms of the local jump components $ \Delta t_{l},  \Delta t_{n}$. The latter components can be further expressed in terms of $\sigma^{S},\, \omega^{S}$ and their derivatives using Eq.~(\ref{Eq:2}). Finally, after all those developments, the boundary integral equations (BIEs) for $\sigma^{S},\, \omega^{S}$ can be obtained from Eqs.~(\ref{Eq:3})-(\ref{Eq:4}) assuming that $\mathbf{x} = \mathbf{y}_0$, with $\mathbf{y}_0 \in L$, see \ref{Section:Appendix A} for details. The resulting system of singular BIEs must be supplemented by the tip conditions of Eq.~(\ref{Eq:5}).

After the system is solved and the values of $\sigma^{S},\, \omega^{S}$ are found, displacements at any point within the domain of interest can be determined using the integral representation for displacements in terms of $\sigma^{S},\, \omega^{S}$. Subsequently, the strains and stresses anywhere within the material system can be computed by proper differentiation of the integral representations for the displacements and using the following relations:
\begin{equation}\label{Eq:8}
    \begin{gathered}
     2\varepsilon_{mn} = u_{m,n} + u_{n,m}\:,\\[10pt]
     \sigma_{mn} = 2\mu\,\varepsilon_{mn} + \frac{(3-\kappa)\mu}{(\kappa-1)}\delta_{mn}\varepsilon_{kk}\,;\quad m,n= 1,2 \,.     
    \end{gathered}
\end{equation}
where $\varepsilon_{mn}$ ($\sigma_{mn}$) are the components of the strain (stress) tensor, $\varepsilon_{kk}$ is the trace of strain tensor, and $u_{m,n} = {\partial u_{m}}/{\partial x_{n}}$.


\section{Numerical Technique}\label{Section:5}

As was emphasized, the final system of singular boundary integral equations involves $\sigma^{S},\,\omega^{S}$ and their derivatives. In addition, the surface may have a shape characterized by varying curvature. These factors render the standard $C^0$ approximation ineffective to solve the system. Therefore, we propose to use IGABEM to address this challenge. In the developed IGABEM algorithm that is described below, we make use of NURBS, see \cite{Cottrell2009}, \cite{PieglLesandTiller1997}, in order to accurately represent the geometry of the surface and approximate the unknown functions $\sigma^{S}$, $\omega^{S}$. Thus, in the following subsection, we provide a brief review of basis functions, such as B-splines, NURBS, and their derivatives.

\subsection{B-Splines and Non-Uniform Rational B-Splines}\label{Sub-Section:5.1}

B-splines are parametric functions of an independent parameter $\xi$ that belongs to the parametric space. The parametric space is defined by the knot vector, which in one dimension is a non-decreasing set of coordinates in the parametric space, $\Xi = [\xi_{1},\xi_{2},\dots,\xi_{n+p+1}]$, where $\xi_{i} \in R$ is the $i^{th}$ knot, $i$ is the knot index, $n$ is the number of basis functions used to construct the B-spline curve and $p$ is the degree of these approximate functions. The knots partition the parametric space into the elements. Element boundaries in physical space are simply images of knot lines under the B-spline mapping. 

The B-spline curve is defined as
\begin{equation}\label{Eq:9}
   \mathbf{C}(\xi) = \sum_{i=1}^{n} N_{i,p}(\xi)\mathbf{P}_{i}\:,
\end{equation}
where $\mathbf{C}(\xi)$ is a vector representing the Cartesian coordinates of a point in physical space, and $\mathbf{P}_{i}$ are the control points, each point is a vector specifying Cartesian coordinates. The B-spline basis functions $N_{i,p}(\xi)$ are defined recursively as follows:\\
for $p=0$
\begin{equation}\label{Eq:10}
   N_{i,0}(\xi) = 
        \begin{cases}
            1 & \text{if }  \xi_{i}\leq \xi < \xi_{i+1} \\
            0 &    \mathrm{otherwise}
        \end{cases}\:,
\end{equation}
for $p=1,2,3,\dots$
\begin{equation}\label{Eq:11}
   N_{i,p}(\xi) = \frac{\xi - \xi_{i}}{\xi_{i+p}-\xi_{i}}N_{i,p-1}(\xi) + \frac{\xi_{i+p+1}-\xi}{\xi_{i+p+1}-\xi_{i+1}} N_{i+1,p-1}(\xi)\:.
\end{equation}

NURBS extend B-splines by incorporating weight functions, thus enabling more accurate approximation of complex geometric shapes. NURBS can be understood from both geometric (see \cite{Farin1999}) and algebraic (see \cite{Cottrell2009}) perspectives. However, in this review, we focus exclusively on the algebraic viewpoint, as it is particularly useful for the development of the algorithm, which is the objective of this study. Accordingly, a NURBS curve is approximated as
\begin{equation}\label{Eq:12}
   \mathbf{C}(\xi) = \sum_{i=1}^{n} R_{i,p}(\xi)\mathbf{P}_{i}\:,
\end{equation}
where $R_{i,p}(\xi)$ denotes the set of NURBS basis functions of degree $p$ at point $\xi$ in parametric space defined as
\begin{equation}\label{Eq:13}
   R_{i,p}(\xi) = \frac{N_{i,p}(\xi)w_{i}}{\sum_{j=1}^{n} N_{j,p}(\xi)w_{j}}\:,
\end{equation}
in which $N_{i,p}$ can be found from  Eqs.~(\ref{Eq:10}) and (\ref{Eq:11}) and $w_{i}$ is the $i$-th weight function. 

From the above expression, the first order derivative of the NURBS basis function can be evaluated as
\begin{equation}\label{Eq:14}
   \frac{\mathrm{d}}{\mathrm{d}\xi}R_{i,p}(\xi) = w_{i}\frac{W(\xi)N_{i,p}^{'}(\xi) - W^{'}(\xi)N_{i,p}(\xi)}{W^{2}(\xi)}\:,
\end{equation}\\
where,
\[
 W(\xi) = \sum_{j=1}^{N} N_{j,p}(\xi)w_{j}\:,\;
 W^{'}(\xi) = \sum_{j=1}^{N} N^{'}_{j,p}(\xi)w_{j}\:,
\]
\[
 \ \; N_{j,p}^{'}(\xi) = \frac{\mathrm{d}}{\mathrm{d}\xi}N_{j,p}(\xi)\:.
\]

The higher-order derivatives of these rational functions can be expressed in terms of lower-order derivatives, e.g., \cite{Cottrell2009}, \cite{Simpson2012}, as
\begin{equation}\label{Eq:15}
   \frac{\mathrm{d}^{k}}{\mathrm{d}\xi^{k}}R_{i,p}(\xi) = \frac{A_{i}^{(k)}(\xi)-\sum_{b=1}^{k} \binom{k}{b} W^{(b)}(\xi)\frac{\mathrm{d}^{(k-b)}}{\mathrm{d}\xi^{(k-b)}}R_{i,p}(\xi)}{W(\xi)}\:,
\end{equation}\\
where, 

\[
\begin{aligned}
A_{i}^{(k)}(\xi) = w_{i}\frac{\mathrm{d}^{k}}{\mathrm{d}\xi^{k}}N_{i,p}(\xi)\:, \\
W^{(k)}(\xi) = \frac{\mathrm{d}^{k}}{\mathrm{d}\xi^{k}}W(\xi)\:, \\
\binom{k}{b} = \frac{k!}{b!\,(k-b)!}\:.
\end{aligned}
\]

\subsection{Isogeometric Boundary Integral Equations}\label{Sub-Section:5.2}

In isogeometric analysis, the unknown fields, the components of the surface stress tensor $\sigma^{S}$ and $\omega^{S}$, are approximated using the same parametric basis functions, NURBS, that define the geometry. 
So, these fields are approximated as
\begin{equation}\label{Eq:16}
    \begin{gathered}
        \sigma^{S} (\xi) = \sum_{i=1}^{n} R_{i,p}(\xi)d_{i} \:, \\[10pt]
        \omega^{S}(\xi) = \sum_{i=1}^{n} R_{i,p}(\xi)q_{i} \:,
    \end{gathered}
\end{equation}
where, $d_{i}$ and $q_{i}$ are the unknown values of $\sigma^{S}$ and $\omega^{S}$ at the control points. 

In the numerical implementation of boundary integral equations (BIEs), the boundary integrals are evaluated over the entire domain by summing the contributions of individual elements. The boundary is parametrized by a NURBS curve, where the parametric domain is partitioned into elements defined by the unique knot values of the knot vector. Thus, each element corresponds to a knot span between two consecutive distinct knots in the parameter space. The NURBS basis functions $R_{i,p}(\xi)$, of degree $p$, exhibit local support; that is, each basis function is nonzero only within a limited parametric interval determined by the knot vector and is identically zero elsewhere. Consequently, at any parametric coordinate $\xi$, only a subset of basis functions -typically $p+1$ - are nonzero and contribute to the approximation. For each element, the subset of nonzero basis functions is identified via a connectivity mapping, which relates the local basis function index $l$ in element $e$ to the corresponding global basis function index $i$. Formally, the local basis functions $N_{l}^{e}(\xi)$ on element $e$ are related to the global NURBS basis functions as $N_{l}^{e}(\xi) = R_{i,p}(\xi)$ with $i = \texttt{conn}(e,l)$, where ``$\texttt{conn}$''
 denotes the connectivity function. This connectivity function encodes the association between the elements and global basis functions, enabling the correct assembly of element contributions in the global system.

Using the above definition of local basis functions, it is now possible to formulate the isogeometric approximations for the geometry and  $\sigma^{S}$, $\omega^{S}$ as 
\begin{equation}\label{Eq:17}
    \begin{gathered}
        \mathbf{C}(\xi) = \sum_{l=1}^{p+1} N_{l}^{e} 
        (\xi)\mathbf{P}_{l}^{e} \:,\\[10pt]
        \sigma^{S} (\xi) = \sum_{l=1}^{p+1} N_{l}^{e}(\xi)d_{l}^{e} \:, \\[10pt]
        \omega^{S}(\xi) = \sum_{l=1}^{p+1} N_{l}^{e}(\xi)q_{l}^{e} \:.
    \end{gathered}
\end{equation}

This parametrization allows the coupled BIEs to be rewritten in parametric form. The BIEs are then evaluated at collocation points in the parametric space, defined by the Greville abscissae \cite{Simpson2013}:
\begin{equation}\label{Eq:18}
\xi_{a}^{'} = \frac{\xi_{a+1} + \xi_{a+2} + \dots + \xi_{a+p}}{p}\:,
\end{equation}
where $a=1, 2, \dots, n$ (as the curve is open-ended) and $\xi_{i}$ denotes the $i^{th}$ knot in the knot vector. 

The parametric form of BIEs for $\sigma^{S}$ and $\omega^{S}$ at the collocation point $\xi_{a}^{'} \in \bar{e}$ can be written as follows:
\begin{multline}\label{Eq:19}
\frac{2\pi\mu (1+\kappa)}{(\lambda_{S} + 2\mu_{S})}
\sum_{l=1}^{p+1} N_{l}^{\bar{e}}(\xi_{a}^{'})d_{l}^{\bar{e}} =\frac{2\pi\mu (1+\kappa)}{(\lambda_{S} + 2\mu_{S})} \sigma_{0}(\xi_{a}^{'}) +\frac{2\pi\mu (1+\kappa)}{(\lambda_{S} + 2\mu_{S})} \Sigma_{3}(\xi_{a}^{'}) \\
+ \sum_{e=1}^{N_{e}} \sum_{l=1}^{p+1} d_{l}^{e} \Bigg\{ 
  \kappa \int_{\Gamma_{e}} \left[
    \cos{\beta^{'}} f_{l}^{e}(\xi)
    + \sin{\beta^{'}} g_{l}^{e}(\xi)
  \right] 
  \frac{
    ( r_{2}\cos{\beta^{'}} - r_{1}\sin{\beta^{'}} )
  }{
    r^2(\xi, \xi_{a}^{'})
  } \mathrm{d}\xi \\
  - 2 \int_{\Gamma_{e}} \left[
    \cos{\beta^{'}} f_{l}^{e}(\xi)
    -\sin{\beta^{'}} g_{l}^{e}(\xi)
  \right] 
   \frac{
    r_{1} r_{2} \left( r_{1} \cos{\beta^{'}} + r_{2}\sin{\beta^{'}}  \right)
  }{
    r^4(\xi, \xi_{a}^{'})
  } \mathrm{d}\xi \\
  + \int_{\Gamma_{e}} \left[
    \sin{\beta^{'}} f_{l}^{e}(\xi)
    + \cos{\beta^{'}}g_{l}^{e}(\xi)
  \right] 
  \frac{
    (r_{1}^2 - r_{2}^2) \left( r_{1} \cos{\beta^{'}} + r_{2}\sin{\beta^{'}}  \right)
  }{
    r^4(\xi, \xi_{a}^{'})
  } \mathrm{d}\xi
  \Bigg\} \\
+ \sigma_{0}(\xi_{a}^{'})\sum_{e=1}^{N_{e}} \sum_{l=1}^{p+1} q_{l}^{e} \Bigg\{
  -\kappa \int_{\Gamma_{e}} \left[
   \sin{\beta^{'}} f_{l}^{e}(\xi)
    - \cos{\beta^{'}} g_{l}^{e}(\xi)
  \right] 
   \frac{
    (r_{2} \cos{\beta^{'}} - r_{1}\sin{\beta^{'}})
  }{
    r^2(\xi, \xi_{a}^{'})
  } \mathrm{d}\xi \\
  - 2 \int_{\Gamma_{e}} \left[
    \sin{\beta^{'}} f_{l}^{e}(\xi)
    + \cos{\beta^{'}}g_{l}^{e}(\xi)
  \right] 
   \frac{
    r_{1} r_{2} \left(r_{1} \cos{\beta^{'}} + r_{2}\sin{\beta^{'}} \right)
  }{
    r^4(\xi, \xi_{a}^{'})
  } \mathrm{d}\xi\\
  - \int_{\Gamma_{e}} \left[
    \cos{\beta^{'}} f_{l}^{e}(\xi)
     -\sin{\beta^{'}}g_{l}^{e}(\xi)
  \right] 
   \frac{
    (r_{1}^2 - r_{2}^2) \left( r_{1} \cos{\beta^{'}} + r_{2}\sin{\beta^{'}} \right)
  }{
    r^4(\xi, \xi_{a}^{'})
  } \mathrm{d}\xi
  \Bigg\}\:,
\end{multline}

\begin{multline}\label{Eq:20}
2\pi\mu(1+\kappa)\sum_{l=1}^{p+1} N_{l}^{\bar{e}}(\xi_{a}^{'})q_{l}^{\bar{e}} = 2\pi\mu(1+\kappa)\Sigma_{4}(\xi_{a}^{'}) \\ 
+  \sum_{e=1}^{N_{e}} \sum_{l=1}^{p+1} d_{l}^{e} \Bigg\{
  -\kappa \int_{\Gamma_{e}} \left[
    -\sin{\beta^{'}} f_{l}^{e}(\xi)
    + \cos{\beta^{'}} g_{l}^{e}(\xi)
  \right] 
\frac{
    ( r_{2} \cos{\beta^{'}} - r_{1}\sin{\beta^{'}} )
  }{
    r^2(\xi, \xi_{a}^{'})
  } \mathrm{d}\xi\\
  - 2 \int_{\Gamma_{e}} \left[
    \sin{\beta^{'}} f_{l}^{e}(\xi)
    + \cos{\beta^{'}} g_{l}^{e}(\xi)
  \right] 
   \frac{
    r_{1} r_{2} \left(  r_{1} \cos{\beta^{'}} + r_{2} \sin{\beta^{'}}  \right)
  }{
    r^4(\xi, \xi_{a}^{'})
  } \mathrm{d}\xi \\
  - \int_{\Gamma_{e}} \left[
    \cos{\beta^{'}} f_{l}^{e}(\xi)
     -\sin{\beta^{'}} g_{l}^{e}(\xi)
  \right] 
   \frac{
    (r_{1}^2 - r_{2}^2) \left(   r_{1} \cos{\beta^{'}} + r_{2} \sin{\beta^{'}} \right)
  }{
    r^4(\xi, \xi_{a}^{'})
  } \mathrm{d}\xi
\Bigg\} \\
+ \sigma_{0} \sum_{e=1}^{N_{e}} \sum_{l=1}^{p+1} q_{l}^{e} \Bigg\{
  \kappa \int_{\Gamma_{e}} \left[
    \cos{\beta^{'}} f_{l}^{e}(\xi)
    + \sin{\beta^{'}} g_{l}^{e}(\xi)
  \right] 
   \frac{
    ( r_{2} \cos{\beta^{'}} - r_{1}\sin{\beta^{'}} )
  }{
    r^2(\xi, \xi_{a}^{'})
  } \mathrm{d}\xi \\
  + 2 \int_{\Gamma_{e}} \left[
    \cos{\beta^{'}} f_{l}^{e}(\xi)
     -\sin{\beta^{'}} g_{l}^{e}(\xi)
  \right] 
  \frac{
    r_{1} r_{2} \left(  r_{1} \cos{\beta^{'}} + r_{2} \sin{\beta^{'}}  \right)
  }{
    r^4(\xi, \xi_{a}^{'})
  } \mathrm{d}\xi \\
  - \int_{\Gamma_{e}} \left[
    \sin{\beta^{'}} f_{l}^{e}(\xi)
    + \cos{\beta^{'}} g_{l}^{e}(\xi)
  \right] 
   \frac{
    (r_{1}^2 - r_{2}^2) \left(  r_{1} \cos{\beta^{'}} + r_{2} \sin{\beta^{'}}  \right)
  }{
    r^4(\xi, \xi_{a}^{'})
  } \mathrm{d}\xi
\Bigg\}\:,
\end{multline}

where
\begin{equation}\label{Eq:21}
   f_{l}^{e}(\xi) = \cos{\beta^{'}}\frac{\mathrm{d}}{\mathrm{d}\xi}N_{l}^{e}(\xi) - \sin{\beta^{'}} \frac{J_{1}(\xi_{a}^{'})}{R(\xi)}N_{l}^{e}(\xi)\:, 
\end{equation}

\begin{equation}\label{Eq:22}
   g_{l}^{e}(\xi) = \sin{\beta^{'}}\frac{\mathrm{d}}{\mathrm{d}\xi}N_{l}^{e}(\xi) + \cos{\beta^{'}}\frac{J_{1}(\xi_{a}^{'})}{R(\xi)}N_{l}^{e}(\xi)\:,
\end{equation}

\begin{equation}\label{Eq:23}
   \Sigma_{3}(\xi_{a}^{'}) = \frac{(\lambda_{S} + 2\mu_{S})}{J_{1}(\xi_{a}^{'})} \left[ -\sin{\beta^{'}} \frac{\mathrm{d}u_{1}^{\infty}(\xi)}{\mathrm{d}\xi}\middle|_{\xi_{a}^{'}}  + \cos{\beta^{'}} \frac{\mathrm{d}u_{2}^{\infty}(\xi)}{\mathrm{d}\xi}\middle|_{\xi_{a}^{'}} \right]\,,
\end{equation}

\begin{equation}\label{Eq:24}
   \Sigma_{4}(\xi_{a}^{'}) = \frac{1}{J_{1}(\xi_{a}^{'})} \left[ \cos{\beta^{'}} \frac{\mathrm{d}u_{1}^{\infty}(\xi)}{\mathrm{d}\xi}\middle|_{\xi_{a}^{'}}  + \sin{\beta^{'}} \frac{\mathrm{d}u_{2}^{\infty}(\xi)}{\mathrm{d}\xi}\middle|_{\xi_{a}^{'}} \right]\:,
\end{equation}

\begin{equation}\label{Eq:25}
   J_{1}(\xi) = \sqrt{\left[\mathbf{C}_{x_{1},\xi}(\xi)\right]^{2} + \left[\mathbf{C}_{x_{2},\xi}(\xi)\right]^{2} }\,.
\end{equation}
In the above set of equations, $\Gamma_{e}$ is the element of integration in the parametric space, $\beta'$ represents the angle formed by the normal at the physical point corresponding to collocation point with the $x$-axis, the terms $\Sigma_{3}(\xi_{a}')$ and $\Sigma_{4}(\xi_{a}')$ correspond to the far-field loading. The $\mathbf{C}_{x_{1}}(\xi)$ and $\mathbf{C}_{x_{2}}(\xi)$ represent the physical coordinates $x_{1}$ and $x_{2}$, respectively, of the point $\xi$ in the parametric domain. Their parametric derivatives, $\mathbf{C}_{x_{1},\xi}(\xi) = {\mathrm{d}}\mathbf{C}_{x_{1}}(\xi)/{\mathrm{d}\xi}$ and $\mathbf{C}_{x_{2},\xi}(\xi) = {\mathrm{d}}\mathbf{C}_{x_{2}}(\xi)/{\mathrm{d}\xi}$, describe the rate of change of the physical coordinates with respect to the parameter $\xi$. At last, $J_{1}(\xi)$ is the Jacobian of transformation from physical space to parametric space.

Additionally, in the above expressions, 
$r_{1} = \mathbf{C}_{x_{1}}(\xi) - \mathbf{C}_{x_{1}}(\xi_{a}^{'})$ and 
$r_{2} = \mathbf{C}_{x_{2}}(\xi) - \mathbf{C}_{x_{2}}(\xi_{a}^{'})$ 
denote the Cartesian components of the relative position vector between the points on the surface $L$ in the physical space corresponding to the parametric coordinates $\xi$ and $\xi_{a}^{'}$; the latter is to the collocation point in parametric space. 
The variable $r = \sqrt{r_{1}^{2} + r_{2}^{2}}$ represents the Euclidean distance between these two points.

Thus, by enforcing the system of BIEs at $n$ collocation points within the parametric domain and substituting the NURBS-based approximations of $\sigma^{S}(\xi)$ and $\omega^{S}(\xi)$, the integral equations are transformed into a set of algebraic equations that involve the unknown coefficients at the control points. The element-wise integrals of kernel functions are numerically evaluated as explained below and assembled into the global matrix $\mathbf{A}$ using the connectivity mapping that relates the local element basis functions to the global control points. Performing such a procedure for all collocation points results in the linear system of the form 
\begin{equation}\label{Eq:26}
\mathbf{A}\mathbf{X} = \mathbf{B}\:,
\end{equation}
where $\mathbf{A}$ is the assembled matrix of the system, $\mathbf{X}$ is the vector of unknown NURBS coefficients representing $\sigma^{S}$ and $\omega^{S}$ at the control points, and $\mathbf{B}$ is the known vector that accounts for the applied external loads and the prescribed surface tension on the surface.

\subsection{Evaluation of the integrals involved}\label{Sub-Section:5.3}
In IGABEM, accurate evaluation of boundary integrals involving kernel functions over element domains is crucial. Depending on the position of the collocation point with respect to the integration element, the integrands may be regular or singular. In this study, we encounter two types of integrals:
(i) regular integrals, when a collocation point is located outside an integration element, and (ii) singular integrals of order $\mathcal{O}\left(\frac{1}{r}\right)$, when a collocation point is located within an integration element.

\begin{figure}[htb]
 	\centering
 	\includegraphics[width=0.8\textwidth]{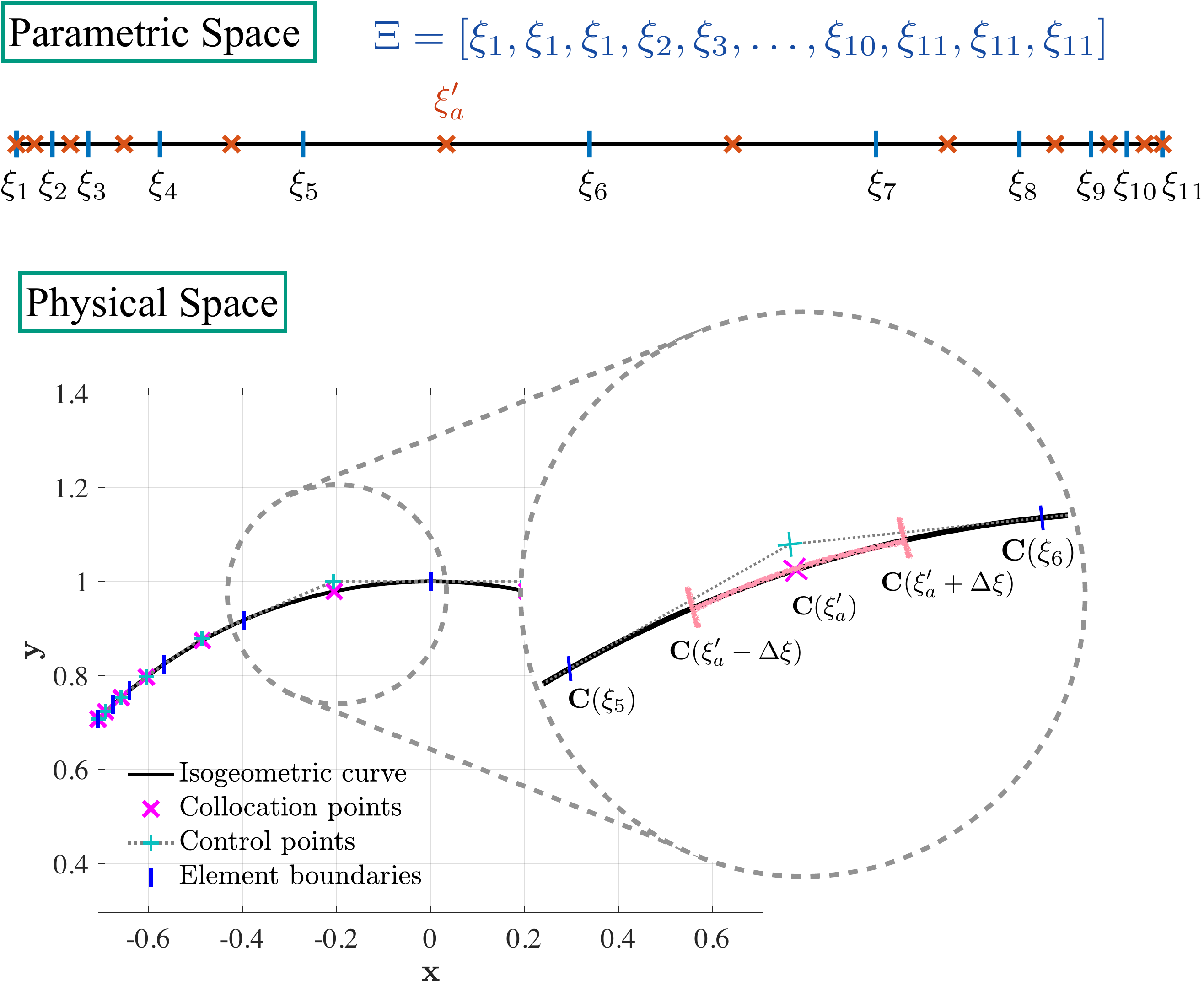}
 	\caption{Illustration of the procedure for evaluation of singular integrals.} 
     \label{Fig:2}
\end{figure}

Using the example of NURBS-defined circular arc of Fig.~\ref{Fig:2}, the strategy for the computation of singular integrals is illustrated. In this example, the coordinates of the arc in physical space, $\mathbf{C}(\xi)$, are mapped from parametric space defined by the knot vector $\Xi$ by using the NURBS basis functions and control points with appropriate corresponding weights. As the BIEs are evaluated at the collocation points $\xi_{a}^{\prime}$ in the parametric space, all integrals of Eqs.~(\ref{Eq:19})-(\ref{Eq:20}) are regular or singular, based on their position with respect to the integration element $\Gamma_{e}$ as mentioned above. 

When $\xi_{a}^{\prime} \notin \Gamma_{e}$, all integrals are regular and can be evaluated using standard Gaussian quadrature rule by appropriately mapping the element $\Gamma_{e}$ in the parent domain of $\hat{\xi}\in[-1,1]$ with the constant Jacobian of transformation
\begin{equation}\label{Eq:Jacobian-2}
J_{2} = \mathrm{d}\xi/\mathrm{d}\hat{\xi}\,.
\end{equation}
When $\xi_{a}^{\prime} \in \Gamma_{e}$, all integrals are singular and special care is needed. In the latter case, $\Gamma_{e}$ is split into sub-intervals, i.e., when $\xi_{a}^{\prime} \in [\xi_{5}, \xi_{6}]$ then
\begin{equation}\label{Eq:27}
\Gamma_{e} = [\xi_{5}, \xi_{a}^{\prime}-\Delta\xi) \cup [\xi_{a}^{\prime}-\Delta\xi, \xi_{a}^{\prime} + \Delta\xi] \cup (\xi_{a}^{\prime}+\Delta\xi,\xi_{6}]\,.
\end{equation}
In the sub-intervals, which do not contain $\xi_{a}^{\prime}$, all integrals are treated as regular ones and evaluated using standard Gaussian quadrature rule. In the remaining sub-interval, which contains $\xi_{a}^{\prime}$, the vector $\mathbf{r} = \mathbf{C}(\xi)- \mathbf{C}({\xi_{a}^{\prime}})$ such that $\xi, \xi_{a}^{\prime} \in  [\xi_{a}^{\prime}-\Delta\xi, \xi_{a}^{\prime} + \Delta\xi] $ is approximated using the Taylor series expansion up to linear term; thus, the distance $r(\xi,\xi_{a}^{\prime})$ can be written as $r(\xi,\xi_{a}^{\prime})= |\xi - \xi_{a}^{\prime}|\, J_{1}(\xi_{a}^{\prime})$. This allows for the reduction of the singular integral to the standard Cauchy Principal Value (C.P.V.) integral of the type 
\[
\int_{\xi_{a}^{\prime}-\Delta\xi}^{\xi_{a}^{\prime}+\Delta\xi} \frac{h(\xi)}{\xi - \xi_{a}^{\prime}}\mathrm{d}\xi\:,
\]
which can further be evaluated using the subtraction of singularity technique (SST), see \cite{Guiggiani1987}, leading to the following representation:
\begin{equation}\label{Eq:28}
\int_{\xi_{a}^{\prime}-\Delta\xi}^{\xi_{a}^{\prime}+\Delta\xi} \frac{h(\xi)}{\xi - \xi_{a}^{\prime}}\mathrm{d}\xi = 
\int_{\xi_{a}^{\prime}-\Delta\xi}^{\xi_{a}^{\prime}+\Delta\xi}
\frac{h(\xi)-h(\xi_{a}^{'})}{\xi - \xi_{a}^{'}} \mathrm{d}\xi \,
+ \mathrm{C.P.V.} \int_{\xi_{a}^{\prime}-\Delta\xi}^{\xi_{a}^{\prime}+\Delta\xi}
\frac{h(\xi_{a}^{'})}{\xi - \xi_{a}^{'}} \mathrm{d}\xi \:,
\end{equation}
where the first term on the right-hand side is evaluated numerically using standard Gaussian quadrature and the second term, the C.P.V. integral, is treated analytically, with the singular behavior explicitly captured by the logarithmic term, refer \cite{Guiggiani1987}.

\subsection{Evaluation of the elastic fields in the material system}\label{Sub-Section:5.4}
After the values of $\sigma^{S},\, \omega^{S}$ are found, all elastic fields at any point outside the surface can be determined using integral representations in terms of $\sigma^{S},\, \omega^{S}$, as explained in Section~\ref{Section:4}. All integrals involved in those representations are regular ones and they can be evaluated using standard quadrature rules.

\section{Numerical Examples}\label{Section:6}

In this section, we present several numerical examples in order to i) validate the proposed IGABEM algorithm using the only two available benchmark solutions and ii) illustrate the influence of curvature variations on the elastic fields. The first benchmark solution involves a surface along a straight segment \cite{Mogilevskaya2021_CST}, while the second solution involves a surface along a circular arc \cite{Han2024_IJSS}. Finally, we present new numerical examples for surfaces of elliptical shapes. 

\subsection{Validation and convergence study }\label{Sub-Section:6.1}

We start validation by considering the problem of a material surface located along a straight segment (a curve characterized by an infinite radius of curvature). The benchmark solutions for the problem are available in \cite{Mogilevskaya2021_CST}, where $\sigma^{S}$ and $\omega^{S}$ were approximated globally using series expansions of second-kind Chebyshev polynomials with square-root weight functions in order to accommodate for the tip conditions of Eq.~(\ref{Eq:5}). In the latter paper, dimensionalization was carried out with respect to the half-length of the segment, and the following dimensionless parameters were introduced: $\widetilde{\sigma}^{S} = \sigma^{S}/\mu a$, $\gamma = 2\mu a/(2\mu_{S}+ \lambda_{S})$, $\widetilde{\sigma}_{0} = \sigma_{0}/\mu a $ and $\widetilde{\sigma}_{ij}= \sigma_{ij}/\mu$, where $a$ denoted the half-length of the straight segment.  

\begin{figure}[htb]
	\centering
	\includegraphics[width=0.75\textwidth]{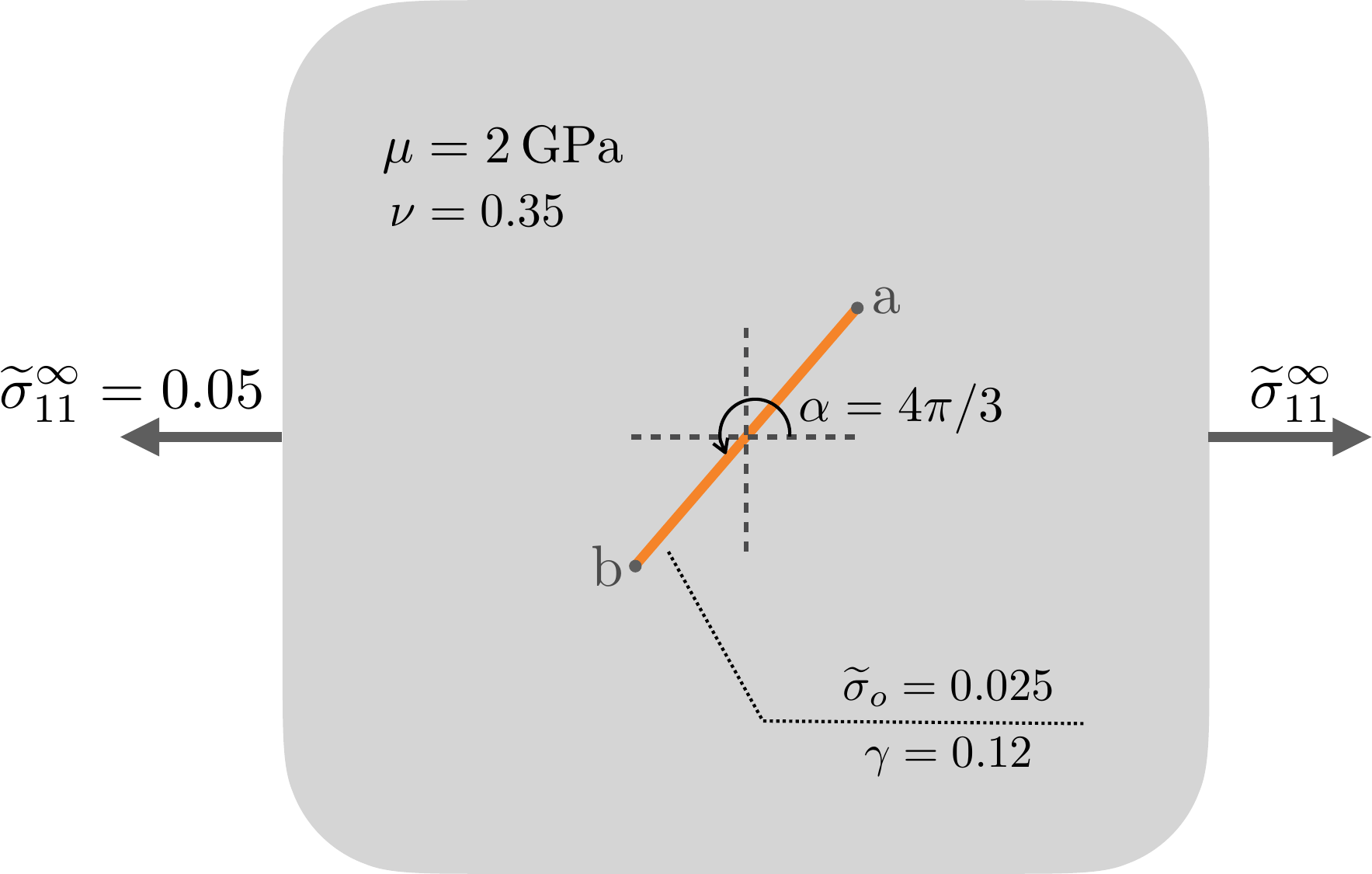}
	\caption{Benchmark problem 1: Surface along the straight segment.}
    \label{Fig:3}
\end{figure}

To enable comparison, we consider a surface represented by the straight segment of length $2a = 10\,\mathrm{nm}$ centered at the origin and located as shown in Fig.~\ref{Fig:3}. The surface is characterized by the dimensionless parameters $\gamma = 0.12$ and $\widetilde{\sigma}_{0} = 0.025$, while the material constants for the matrix are chosen to be $\mu = 2\,\mathrm{GPa}$ and $\nu = 0.35$. The only nonzero component of the far-field load is taken to be $\widetilde{\sigma}_{11}^{\infty} = 0.05$.

To ensure convergence of global approximations, the results reported in \cite{Mogilevskaya2021_CST} were obtained using 40 terms of the Chebychev series expansion for $\widetilde{\sigma}^{S}$ and $\omega^{S}$, and 100 collocation points to assemble the linear algebraic equation system for unknown coefficients of the series.

\begin{figure}[htbp]
  \centering
  \begin{minipage}[b]{0.49\textwidth}
      \centering
      \includegraphics[angle=0, width=1\textwidth]{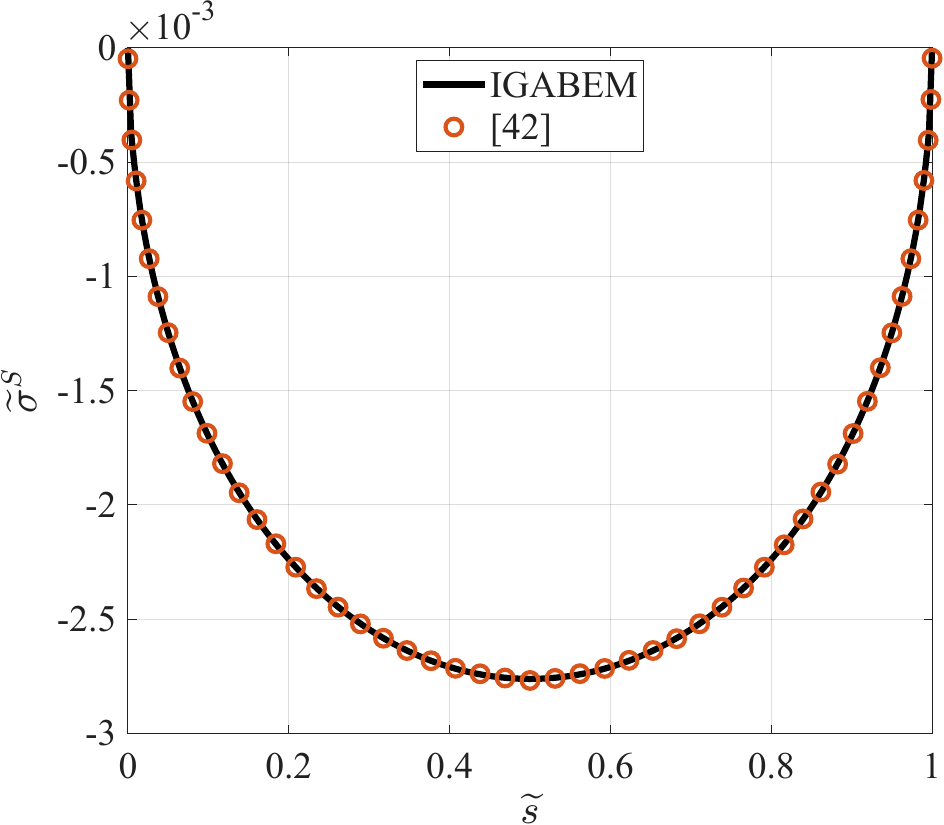}
      \subcaption{}
      \label{Fig:4a}
  \end{minipage}
  \hfill
  \begin{minipage}[b]{0.49\textwidth}
      \centering
      \includegraphics[angle=0, width=1\textwidth]{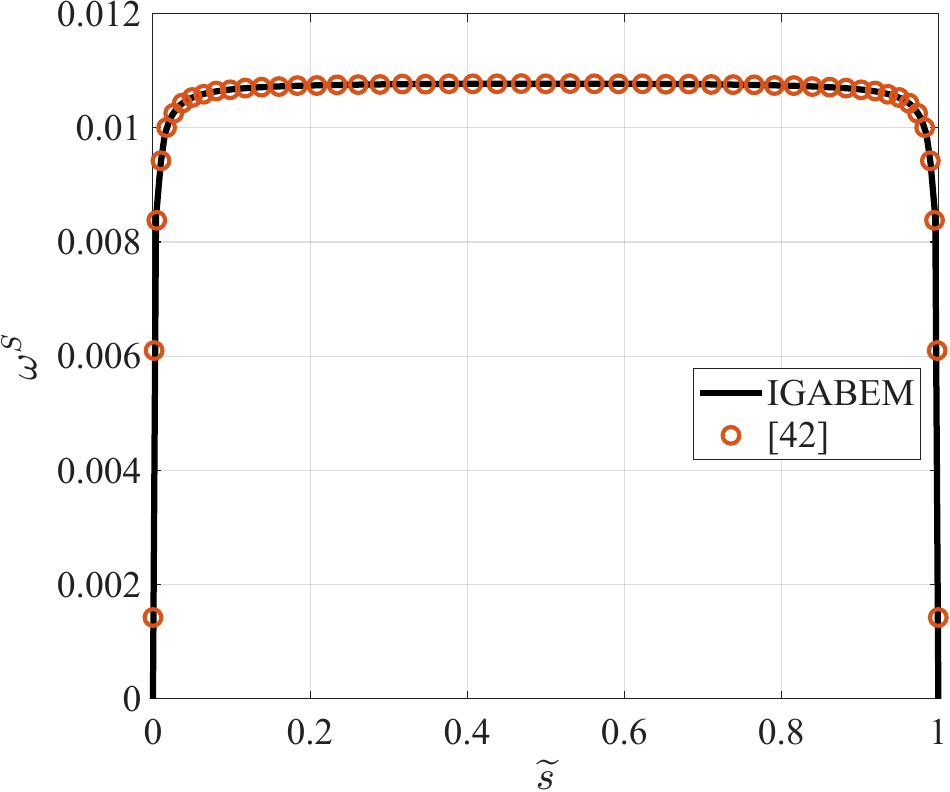}
      \subcaption{}
      \label{Fig:4b}
  \end{minipage}
  \caption{ Comparison of the results for (a) $\widetilde{\sigma}^{S}$ (left)  and (b) $\omega^{S}$ (right) with benchmark solutions for the surface along the inclined straight segment.}
  \label{Fig:4}
\end{figure}

In the IGABEM implementation, the straight segment was modeled using a quadratic (second-degree) NURBS curve. Convergence was achieved with $N_{e}=50$ and 200 Gaussian quadrature points for numerical integration, see \ref{Section:Appendix B}. On Fig.~\ref{Fig:4} the distributions of $\widetilde{\sigma}^{S}$ and $\omega^{S}$ are plotted as functions of the normalized arc length $\widetilde{s} = s/2a$. It can be seen that there is good agreement between the results obtained with the IGABEM algorithm and those reported in \cite{Mogilevskaya2021_CST}. Unlike in the later paper, where the unknowns were approximated globally, the IGABEM algorithm did not require the use of spectral filtering techniques, i.e. \cite{Sarra2006}, to accommodate the tip conditions.

\begin{figure}[H]
  \centering
  \begin{subfigure}[b]{0.49\textwidth}
      \centering
      \includegraphics[angle=0, width=1\textwidth]{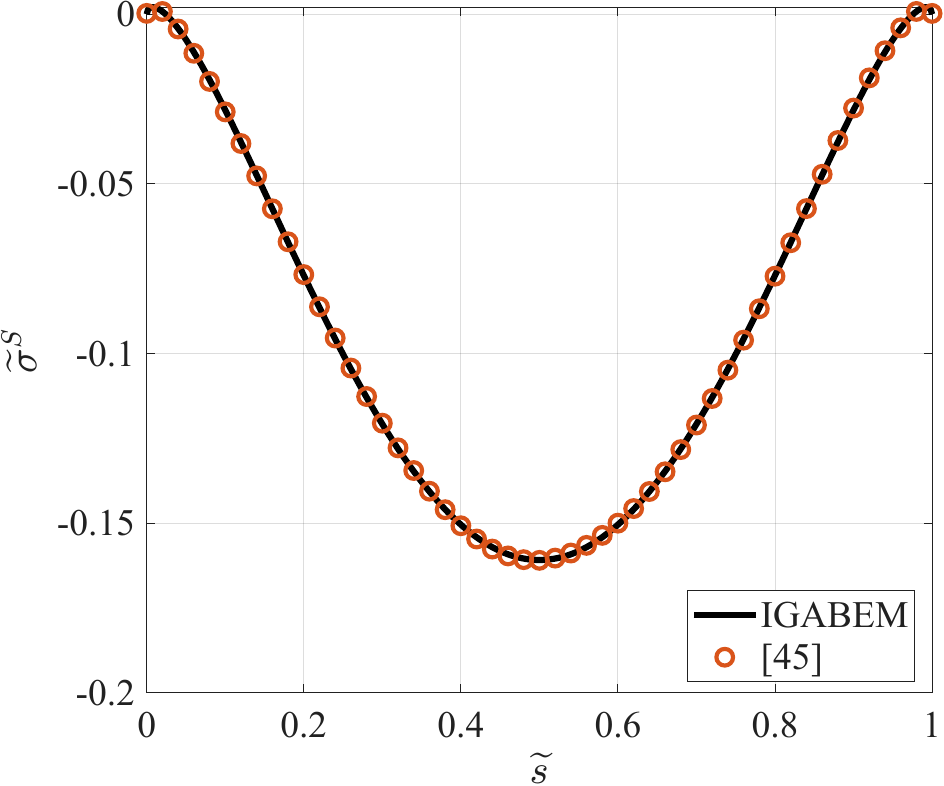}
      \subcaption{}
      \label{Fig:6a}
  \end{subfigure}
  \hfill
  \begin{subfigure}[b]{0.49\textwidth}
      \centering
      \includegraphics[angle=0, width=1\textwidth]{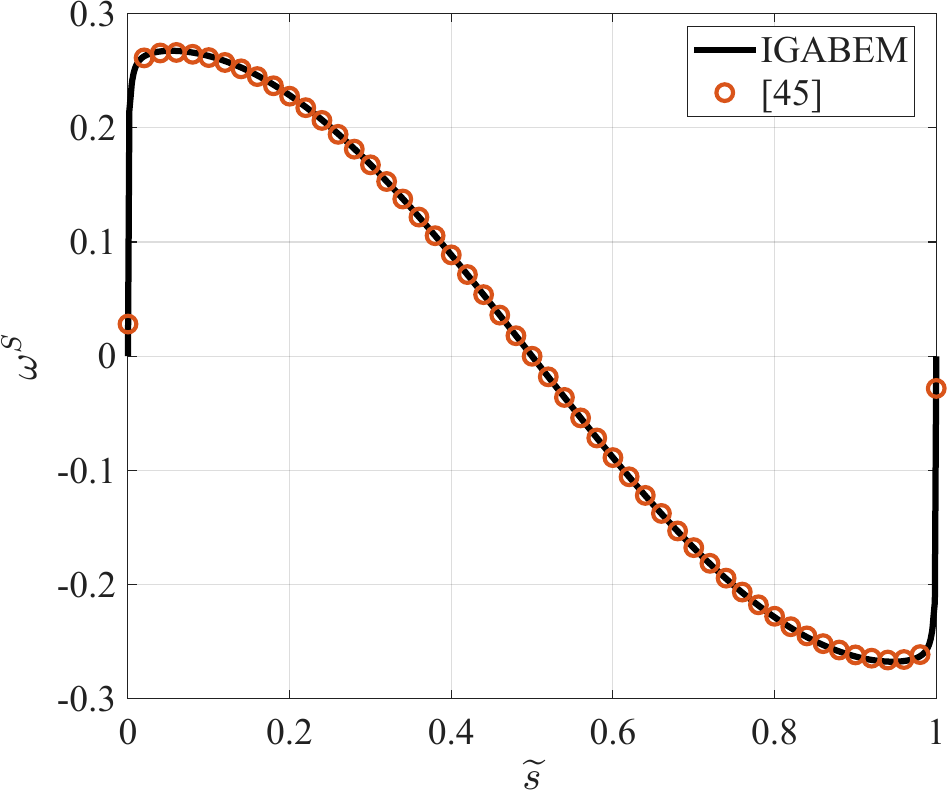}
      \subcaption{}
      \label{Fig:6b}
  \end{subfigure}
  \caption{ Comparison of the results for (a) $\widetilde{\sigma}^{S}$ (left)  and (b) $\omega^{S}$ (right) with benchmark solutions for the surface along the circular arc.}
  \label{Fig:6}
\end{figure}

To further validate our algorithm, we consider the second benchmark problem where the surface has a circular arc shape, i.e., it has a constant but finite radius of curvature. The benchmark solutions for the problem were reported in \cite{Han2024_IJSS}. There, the dimensionalization was performed with respect to the semi-arc length $s = R\theta/2$, and the following dimensionless parameters were introduced: $\gamma = \mu R \theta/(2\mu_{S}+\lambda_{S})$, $\widetilde{\sigma}^{S} = 2\sigma^{S}/\mu R \theta$, $\widetilde{\sigma}_{0} = 2\sigma_{0}/\mu R \theta$, and $\widetilde{\sigma}_{ij}= \sigma_{ij}/\mu$, where $\theta = \beta_{2}-\beta_{1}$ and $\beta$ denotes the angle between the normal at a point on the arc and the $x$-axis. The unknown values of $\widetilde{\sigma}^{S}$ and $\omega^{S}$  were approximated globally using a series of trigonometric functions multiplied by square root weight functions to accommodate the tip conditions of Eq.~(\ref{Eq:5}). To obtain the solution shown in Fig.~\ref{Fig:6}, 40 terms of the truncated trigonometric series and 800 Gaussian quadrature points were considered.
\begin{figure}[htb]
	\centering
	\includegraphics[width=0.525\textwidth]{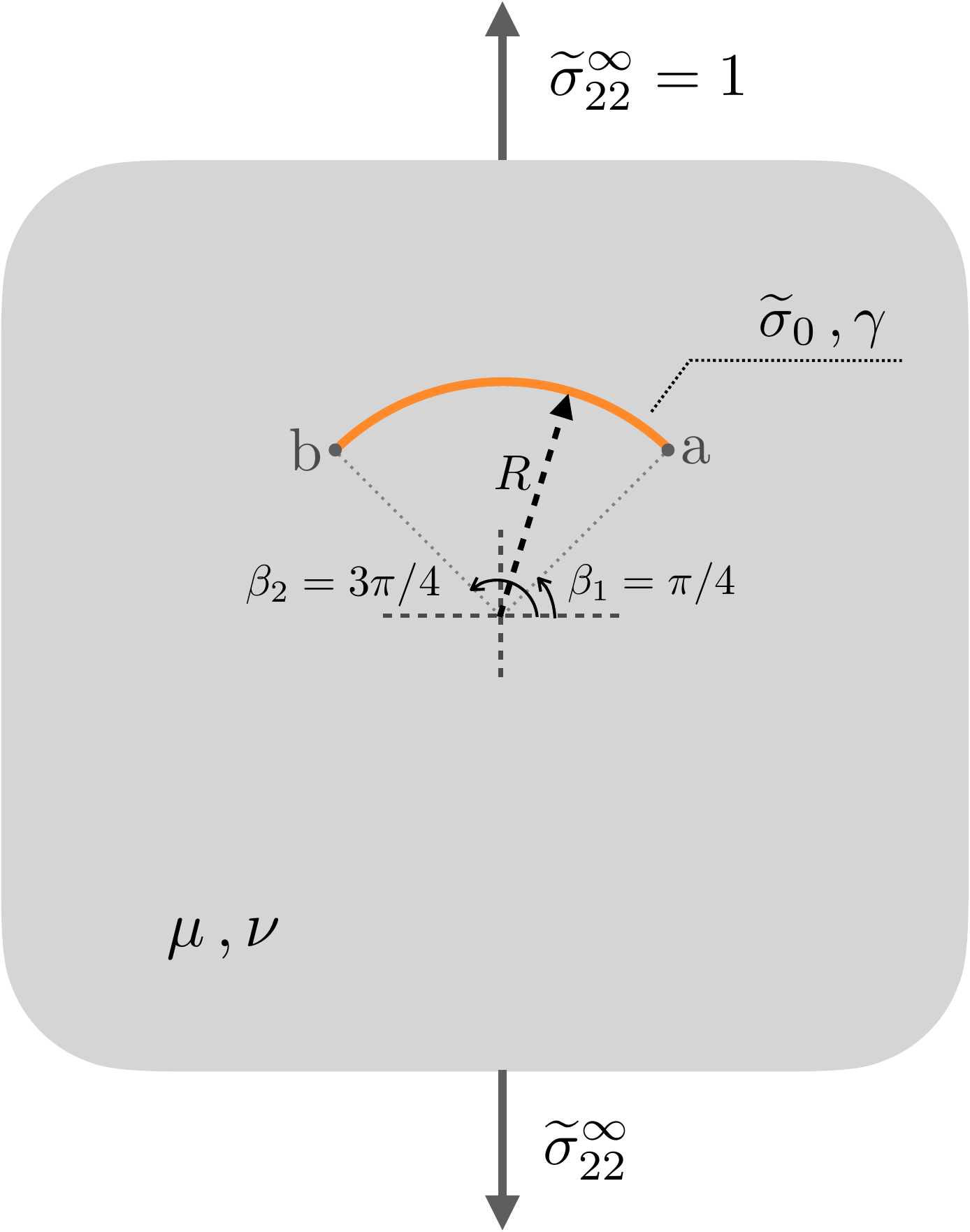}
	\caption{Benchmark problem 2: surface along the circular arc.}
    \label{Fig:5}
\end{figure}

For comparison, we considered a surface along the circular arc with the center at the origin of the Cartesian coordinate system defined by $\beta = [\pi/4,\, 3\pi/4]$, where $\beta$ was the angle between the $x$-axis and the normals at the tips of the surface as shown in Fig.~\ref{Fig:5}. The radius of the arc was taken to be $1\,\mathrm{nm}$. We chose dimensionless surface parameters as $\gamma = 1$ and $\widetilde{\sigma}_{0}=0.01$. The Poisson’s ratio for the matrix was taken to be $\nu = 0.33$, and the dimensionless load was assumed to be $\widetilde{\sigma}_{22} = 1$.

In our simulations, convergence was achieved using a quadratic NURBS approximation of the circular arc, discretized with 50 elements, see \ref{Section:Appendix B}. For the evaluation of regular integrals, 200 Gaussian quadrature points were used. A graded mesh was adopted: the discretization was refined near the arc tips to capture localized effects, while a coarser mesh was used in the midsection of the curve. Fig.~\ref{Fig:6} shows the plots of $\widetilde{\sigma}^{S}$ and $\omega^{S}$ versus the normalized arc-length defined as $\widetilde{s}= s/2\ell$, where $\ell$ is the semi-arc length of the circular arc, i.e, $\ell = \pi/4$. The IGABEM and benchmark results of \cite{Han2024_IJSS} are compared for $\widetilde{\sigma}^{S}$ and $\omega^{S}$ in  Fig.~\ref{Fig:6a} and Fig.~\ref{Fig:6b}, respectively. In both cases, the IGABEM predictions (solid line) are in excellent agreement with the solution reported in \cite{Han2024_IJSS} (symbols), demonstrating the accuracy and reliability of the proposed formulation. The close overlap persists even near the tips, confirming the ability of IGABEM to capture both smooth variations and localized tip effects.

\begin{figure}[H]
  \centering
  \begin{subfigure}[b]{0.48\textwidth}
      \centering
      \includegraphics[angle=0, width=1\textwidth]{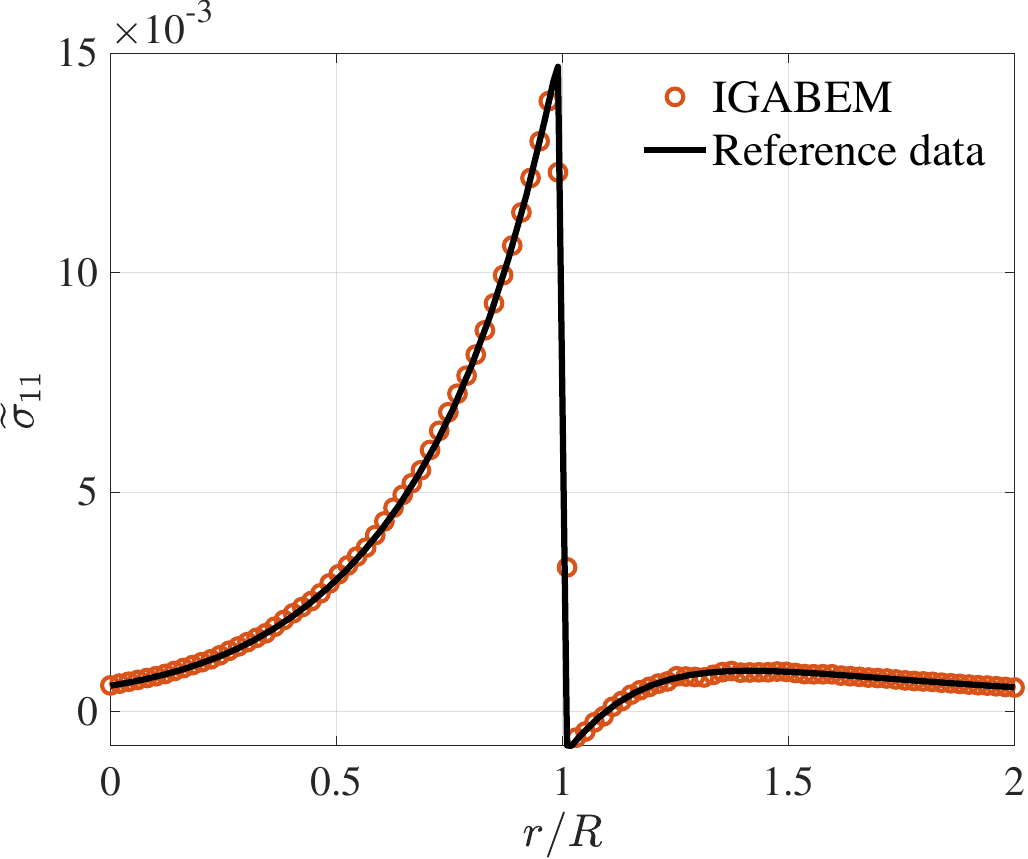}
      \subcaption{}
      \label{Fig:7a}
  \end{subfigure}
  \hfill
  \begin{subfigure}[b]{0.49\textwidth}
      \centering
      \includegraphics[angle=0, width=1\textwidth]{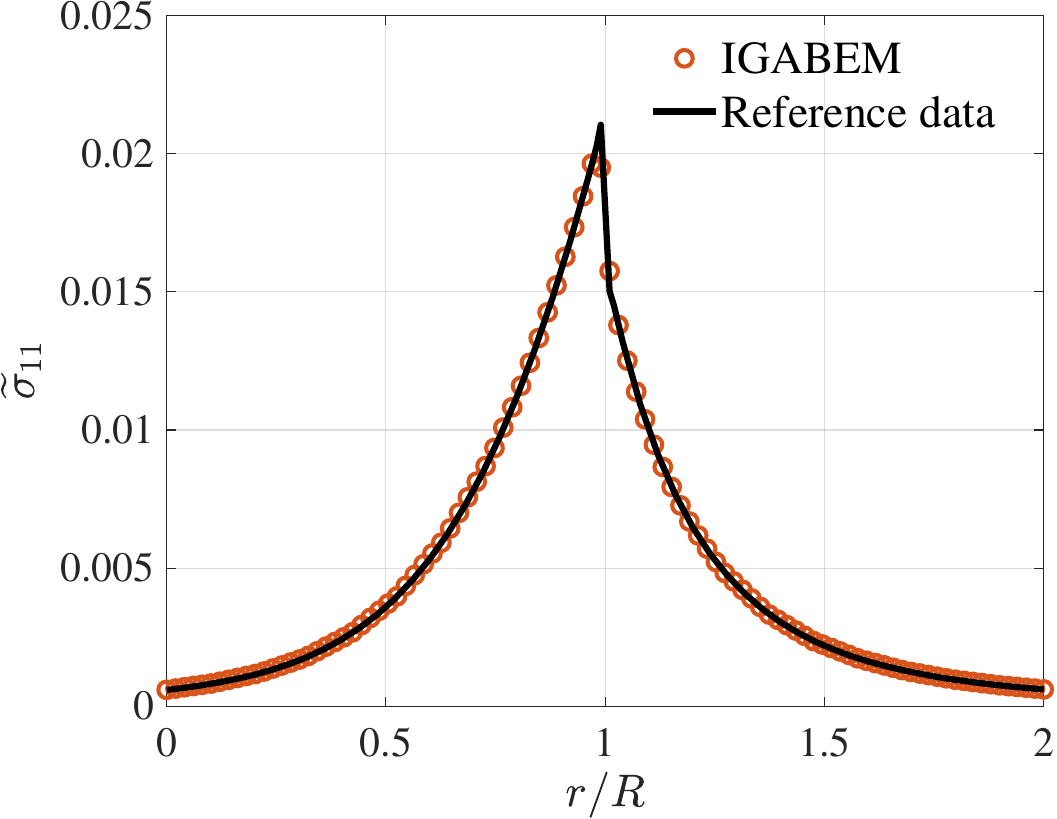}
      \subcaption{}
      \label{Fig:7b}
  \end{subfigure}
  \caption{Comparison of the dimensionless Cauchy stress $\widetilde{\sigma}_{11}$ obtained by IGABEM (symbols) and reference data (solid line) along radial lines at (a) $\beta = 3\pi/8$ (left) and (b) $\beta = \pi/2$ (right).}
  \label{Fig:7}
\end{figure}

\begin{figure}[H]
  \centering
  \begin{subfigure}[b]{0.485\textwidth}
      \centering
      \includegraphics[angle=0, width=1\textwidth]{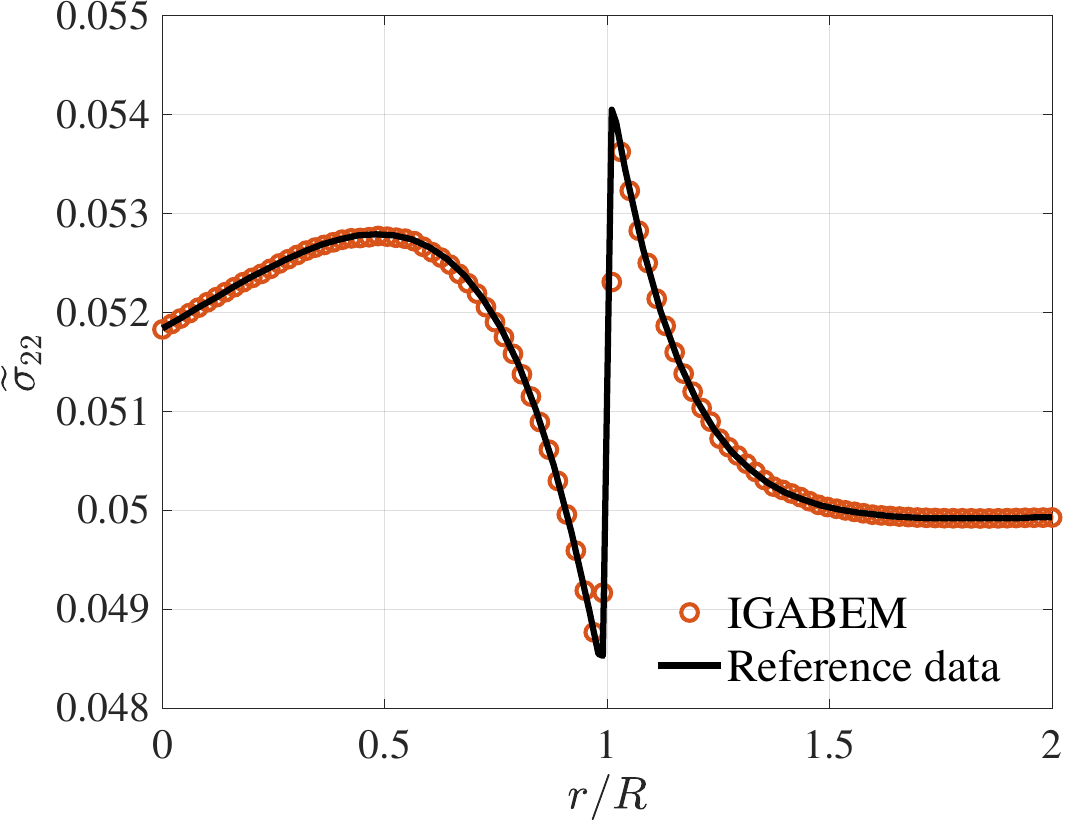}
      \subcaption{}
      \label{Fig:9a}
  \end{subfigure}
  \hfill
  \begin{subfigure}[b]{0.485\textwidth}
      \centering
      \includegraphics[angle=0, width=1\textwidth]{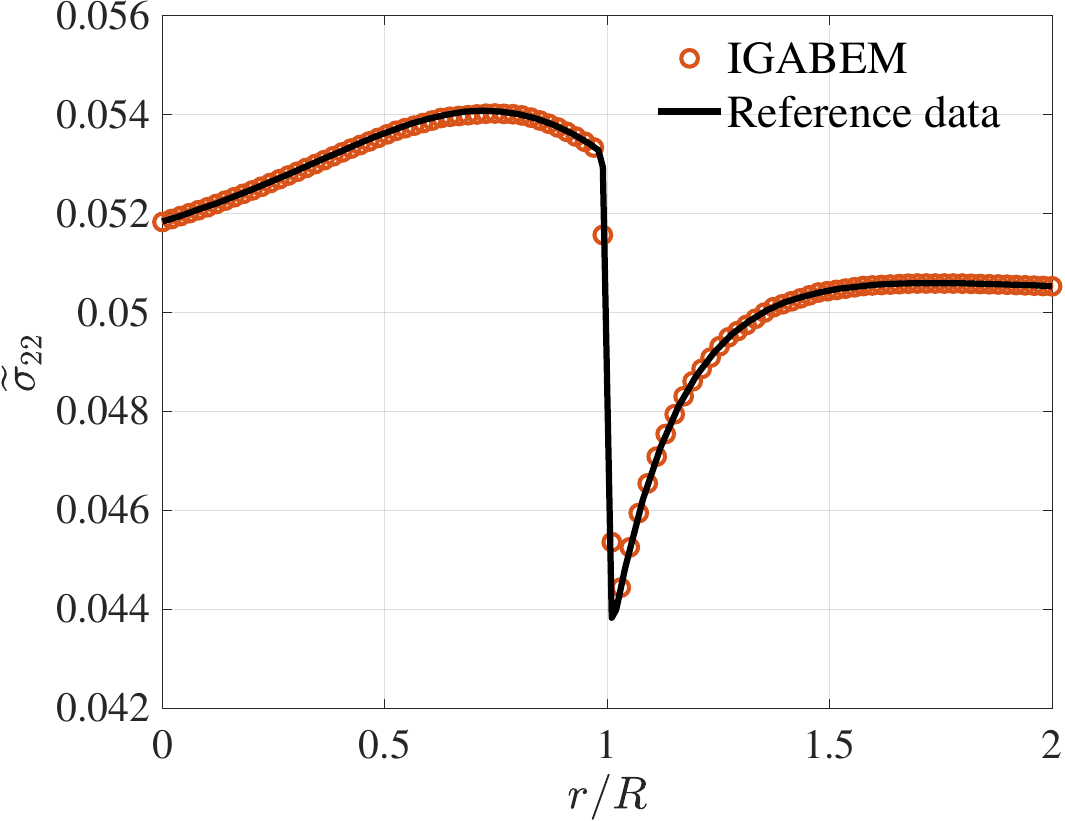}
      \subcaption{}
      \label{Fig:9b}
  \end{subfigure}
  \caption{Comparison of the dimensionless Cauchy stress $\widetilde{\sigma}_{22}$ obtained by IGABEM (symbols) and reference data (solid line) along radial lines at (a) $\beta = 3\pi/8$ (left) and (b) $\beta = \pi/2$ (right). }
  \label{Fig:9}
\end{figure}

Furthermore, Figs.~\ref{Fig:7}\,–\ref{Fig:9} present comparisons of the dimensionless Cauchy stress components, $\widetilde{\sigma}_{ij}$ obtained via IGABEM with the reference data obtained by global approximation-based algorithm of ~\cite{Han2024_IJSS}. The arc with radius $R = 5\,\mathrm{nm}$ is still centered at the origin of the Cartesian coordinate system and defined by $\beta = [\pi/4,\, 3\pi/4]$. The bulk matrix is epoxy with $\mu = 2\,\mathrm{GPa}$ and $\nu = 0.35$. The surface parameters and far-field loading are taken as $\gamma = 0.12$, $\widetilde{\sigma}_{0} = 0.025$, and $\widetilde{\sigma}^{\infty}_{22} = 0.05$.  

\begin{figure}[H]
  \centering
  \begin{subfigure}[b]{0.49\textwidth}
      \centering
      \includegraphics[angle=0, width=1\textwidth]{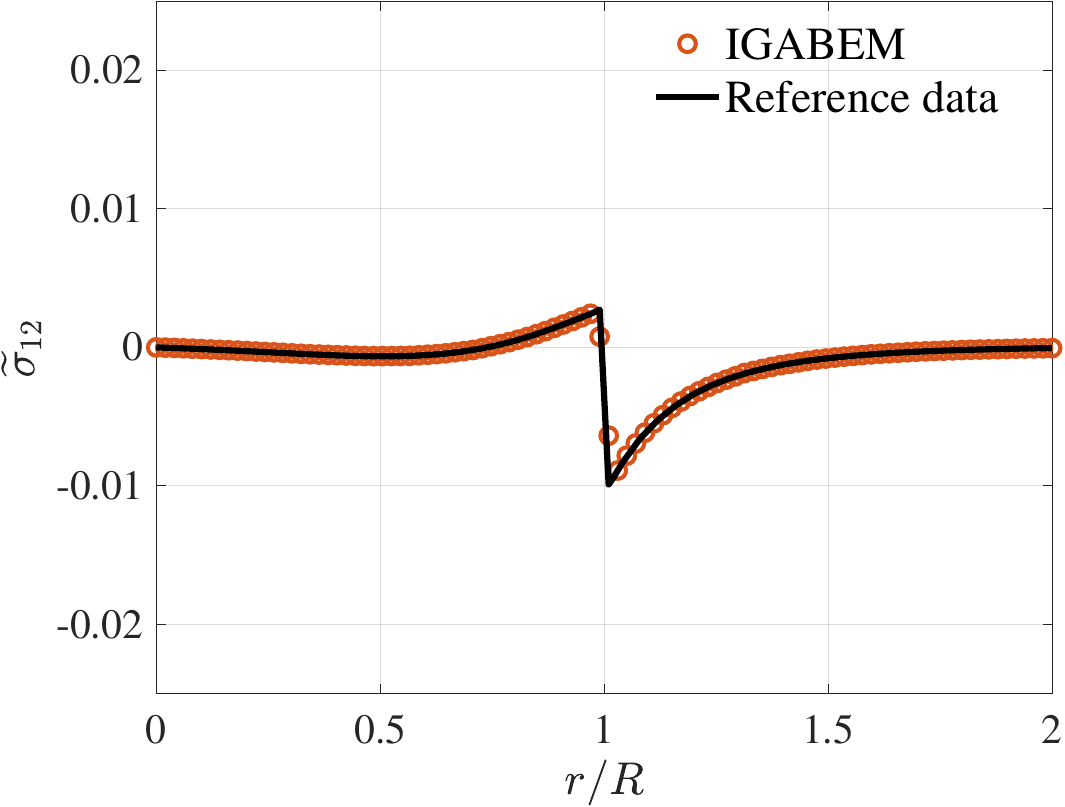}
      \subcaption{}
      \label{Fig:8a}
  \end{subfigure}
  \hfill
  \begin{subfigure}[b]{0.49\textwidth}
      \centering
      \includegraphics[angle=0, width=1\textwidth]{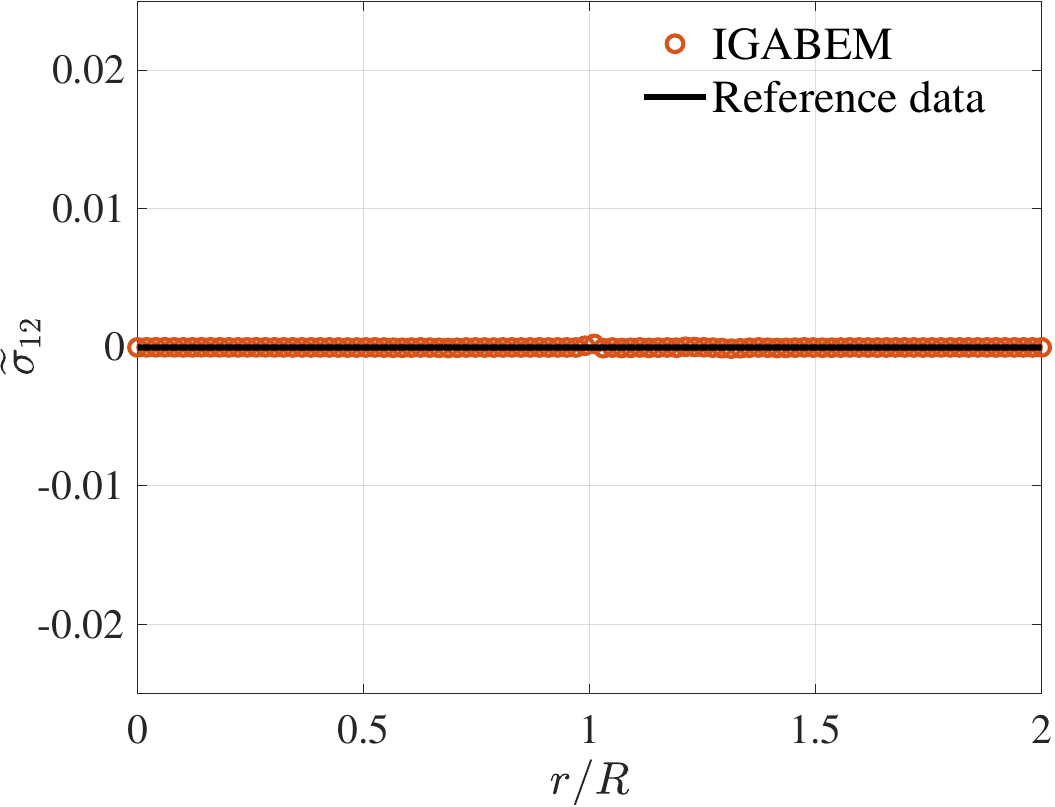}
      \subcaption{}
      \label{Fig:8b}
  \end{subfigure}
  \caption{Comparison of the dimensionless Cauchy stress $\widetilde{\sigma}_{12}$ obtained by IGABEM (symbols) and reference data (solid line) along radial lines at (a) $\beta = 3\pi/8$ (left) and (b) $\beta = \pi/2$ (right). }
  \label{Fig:8}
\end{figure}

In Figs.~\ref{Fig:7}–\ref{Fig:9}, the stress components are plotted along two radial directions, namely $\beta = 3\pi/8$ and $\beta = \pi/2$, as functions of the normalized radial distance $r/R$ measured from the origin of the Cartesian coordinate system.  As can be seen in the figures, good agreement with the reference data is achieved for all stress components.

In addition to the examples discussed above, we have also validated our results using the remaining numerical examples presented in \cite{Mogilevskaya2021_CST} and \cite{Han2024_IJSS}, which involved different material parameters and loading conditions. In all examples, the agreement between the results obtained with the two approaches was also good.

\subsection{Study of curvature effects}\label{Sub-Section:6.2}
\begin{figure}[htb]
	\centering
	\includegraphics[width=0.75\textwidth]{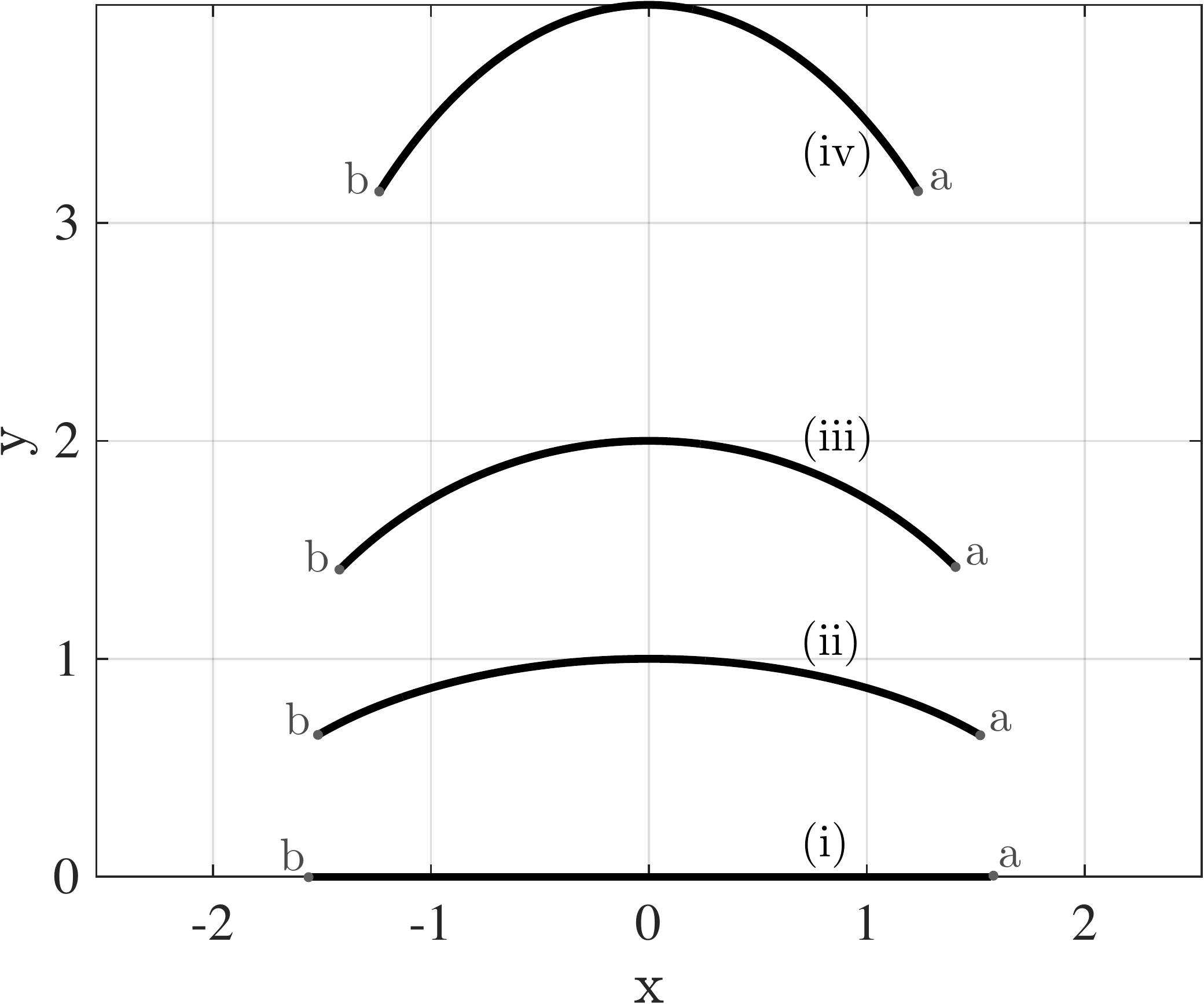}
	\caption{The four surfaces under the study.}
    \label{Fig:10}
\end{figure}

In this study, we consider four representative cases: a surface along $\text{(i)}$ a straight segment, $\text{(ii)}$ an elliptic curve with its major axis oriented horizontally, $\text{(iii)}$ a circular arc, and $\text{(iv)}$ an elliptic curve with its major axis oriented vertically. To perform the study, the length of each surface was assumed to be the same and equal to $2\ell$. In all four cases, the curve on parts of which each surface is located, see Fig.~\ref{Fig:10}, can be described by the following equation:

\begin{equation}\label{Eq:30}
\frac{x^{2}}{a^{2}} + \frac{y^{2}}{b^{2}} = 1\:,
\end{equation}
where $1/b^2=0$ for case $\text{(i)}$; $a>b$ for case $\text{(ii)}$; $a=b$ for case $\text{(iii)}$; and $a<b$ for case $\text{(iv)}$. 

To allow for comparison, all the governing parameters were normalized with respect to the half length of the surface, $\ell$, which, as mentioned, was the same in each case. The resulting non-dimensional parameters are 
\begin{equation}\label{Eq:31}
\gamma = \frac{2\mu\ell}{2\mu_{S} + \lambda_{S}}\,,\, \widetilde{\sigma}^{S}= \frac{\sigma^{S}}{\mu\ell}\,,\,\widetilde{\sigma}_{0} = \frac{\sigma_{0}}{\mu\ell}\,,\,\widetilde{\sigma}^{\infty}_{ij} = \frac{\sigma^{\infty}_{ij}}{\mu}\:.
\end{equation}
For illustration purposes, we consider an epoxy matrix with $\mu = 2\, \mathrm{GPa}$, $\nu = 0.35$. The surface in each case is characterized as follows: $\text{(i)}$ a straight segment of length $2\ell = \pi\,\mathrm{nm}$ defined by the endpoints $\mathbf{a} = (\pi/2,0)$ and $\mathbf{b} = (-\pi/2, 0)$; $\text{(ii)}$ an elliptical curve with $a = 2\,\mathrm{nm}$, $b = 1\,\mathrm{nm}$, with endpoints $\mathbf{a} = (1.5149,\,0.6529)$ and $\mathbf{b} = (-1.5149,\,0.6529)$; $\text{(iii)}$ a circular arc with $a = b = 2\,\mathrm{nm}$, with endpoints $\mathbf{a} = (\sqrt{2},\sqrt{2})$ and $\mathbf{b} = (-\sqrt{2}, \sqrt{2})$; and $\text{(iv)}$ an elliptical curve with $a = 2\,\mathrm{nm}$, $b = 4\,\mathrm{nm}$ with endpoints $\mathbf{a} = (1.2374,\,3.1425)$ and $\mathbf{b} = (-1.2374,\,3.1425)$. The two-dimensional elastic properties of each surface are chosen as: $\gamma = 0.12$ and $\widetilde{\sigma}_{0} = 0.025$. We consider only one non-vanishing component of the far-field load, $\sigma^{\infty}_{22} = 100\,\mathrm{MPa}$, which corresponds to the dimensionless value $\widetilde{\sigma}^{\infty}_{22} = 0.05$.

\begin{figure}[H]
  \centering
  \begin{subfigure}[b]{0.49\textwidth}
      \centering
      \includegraphics[angle=0, width=1\textwidth]{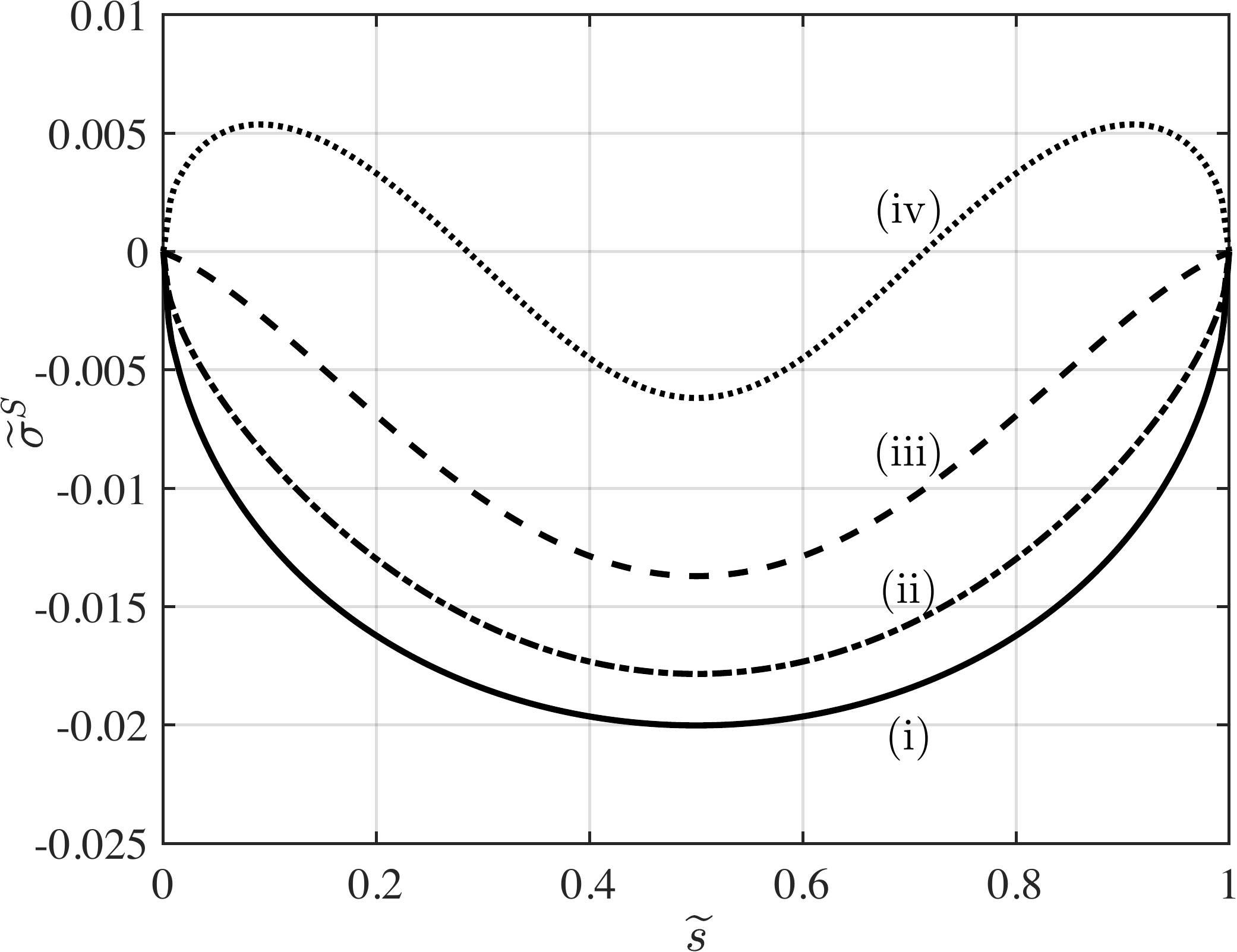}
      \subcaption{}
      \label{Fig:11.a}
  \end{subfigure}
  \hfill
  \begin{subfigure}[b]{0.49\textwidth}
      \centering
      \includegraphics[angle=0, width=1\textwidth]{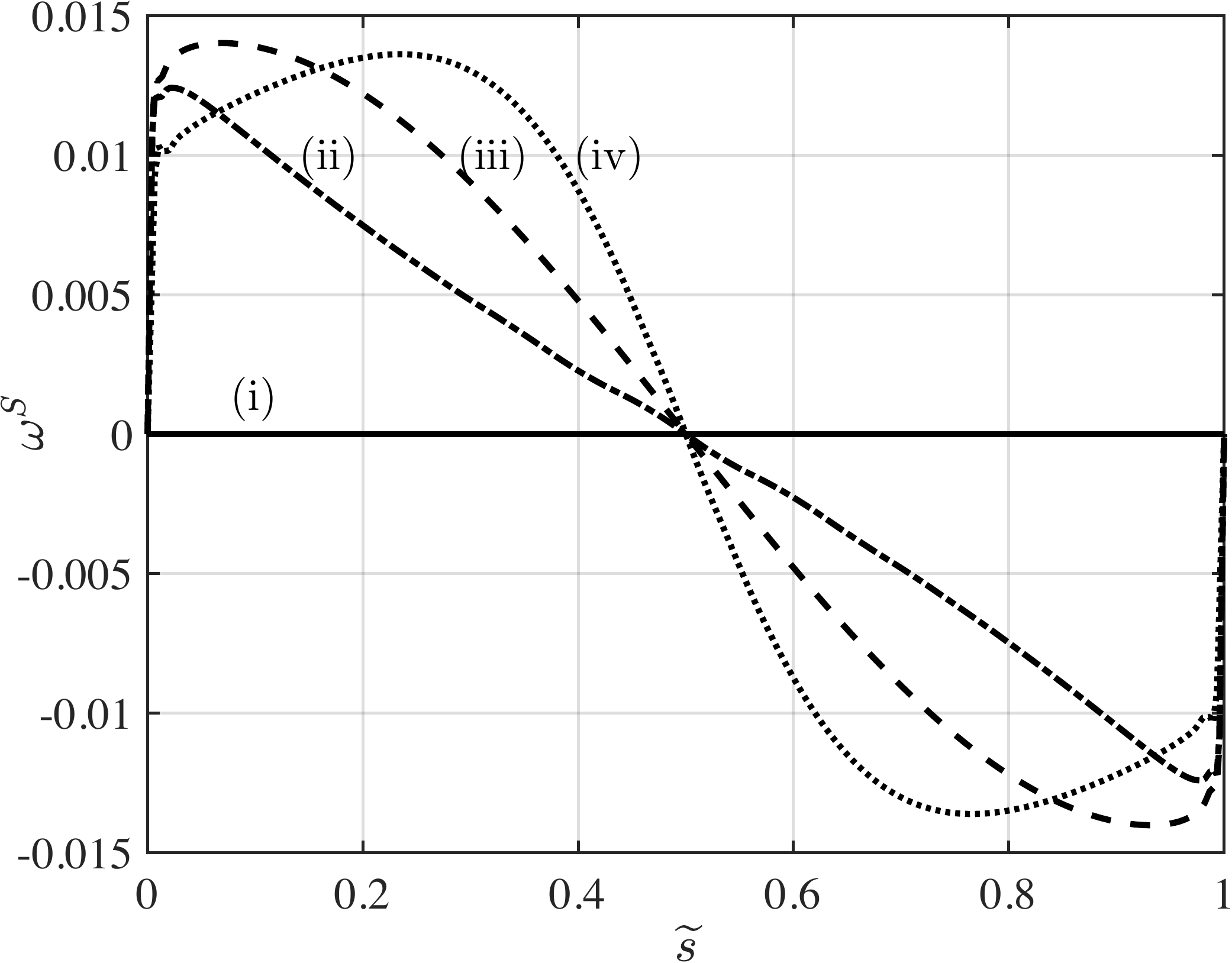}
      \subcaption{}
      \label{Fig:11.b}
  \end{subfigure}
  \caption{ The plots of (a) $\widetilde{\sigma}^{S}$ (left) and (b) $\omega^{S}$ (right) versus normalized arc length for the four surfaces under the study.}
  \label{Fig:11}
\end{figure}

The distributions of $\widetilde{\sigma}^{S}$ and $\omega^{S}$ along different surfaces are plotted in Fig.~\ref{Fig:11} as functions of the normalized arc-length $\widetilde{s} = s/2\ell$. 

As seen in Fig.~\ref{Fig:11.a}, all the plots for $\widetilde{\sigma}^{S}$ are symmetric with respect to $\widetilde{s}=0.5$, as expected. It can be concluded that the stresses $\widetilde{\sigma}^{S}$ in cases $\text{(i)}$-$\text{(iii)}$ are compressive, while gradually converging to zero at the tips. The maximum absolute values of $\widetilde{\sigma}^{S}$ for the three cases are reached at $\widetilde{s} = 0.5$ with the highest absolute values observed for case $\text{(i)}$. The behavior of the plot of $\widetilde{\sigma}^{S}$ for case $\text{(iv)}$ is qualitatively different. The compressive stresses $\widetilde{\sigma}^{S}$ are now localized around $\widetilde{s} = 0.5$ and the distribution rapidly transitions to tensile stresses, approaching zero at the tips more rapidly. Unlike in three former cases, the plot for the latter case has three distinct extrema, one minimum at $\widetilde{s} = 0.5$, and two maxima, roughly around $\widetilde{s} = 0.1$ and $\widetilde{s} = 0.9$. 

The corresponding plots of $\omega^{S}$, Fig.~\ref{Fig:11.b}, for the cases $\text{(ii)}$-$\text{(iv)}$ are antisymmetric with respect to $\widetilde{s}=0.5$, while $\omega^{S}=0$ along the straight surface, case $\text{(i)}$, as expected. The absolute values of $\omega^{S}$ increase as one moves from case $\text{(ii)}$ to case $\text{(iii)}$, abruptly transitioning to zero in the narrow vicinities of the tips, where the non-zero extrema values are reached. The variations in the plot of $\omega^{S}$ for case $\text{(iv)}$ are more pronounced compared to cases $\text{(ii)}$-$\text{(iii)}$ in which the values of $\omega^{S}=0$ varied almost linearly in large parts of the surfaces. For case $\text{(iv)}$, the transition to zero at the tips occurs more gradually and the non-zero extrema values are reached farther away from the tips. Collectively, all of these results emphasize the decisive role of surface curvature variations in distributions of the surface stress components.

To illustrate the effects of curvature variations on the stresses within the matrix, we also present contour plots of the relative dimensionless Von-Mises stress. Here, the relative values are obtained by normalizing with respect to the dimensionless Von-Mises stress corresponding to the case without a surface.  The dimensionless Von-Mises stress under plane strain conditions is defined as
\begin{equation}\label{Eq:32}
    \begin{gathered}
        \widetilde{\sigma}_{v} = \sqrt{\frac{1}{2} \left[ (\widetilde{\sigma}_{11} - \widetilde{\sigma}_{22})^{2} + (\widetilde{\sigma}_{22} - \widetilde{\sigma}_{33})^{2} + (\widetilde{\sigma}_{33} - \widetilde{\sigma}_{11} )^{2} \right] + 3(\widetilde{\sigma}_{12})^{2} } \:, \\[10pt]
    \widetilde{\sigma}_{33} = \nu(\widetilde{\sigma}_{11} + \widetilde{\sigma}_{22})\:.
    \end{gathered}
\end{equation}

\begin{figure}[htbp]
    \centering
    \begin{subfigure}{0.495\textwidth}
        \centering
        \includegraphics[width=\linewidth]{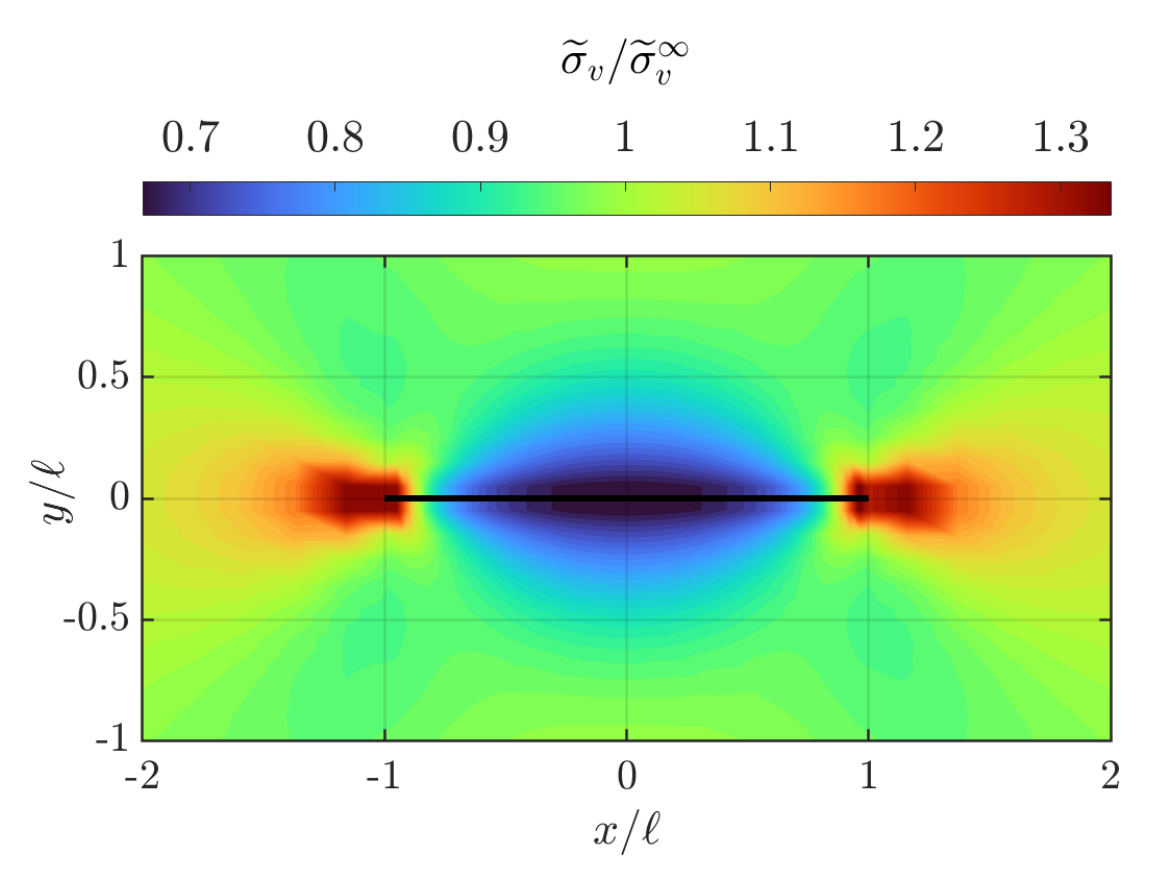}
        \caption{Case (i)}
        \label{Fig:12.a}
    \end{subfigure}
    \begin{subfigure}{0.495\textwidth}
        \centering
        \includegraphics[width=\linewidth]{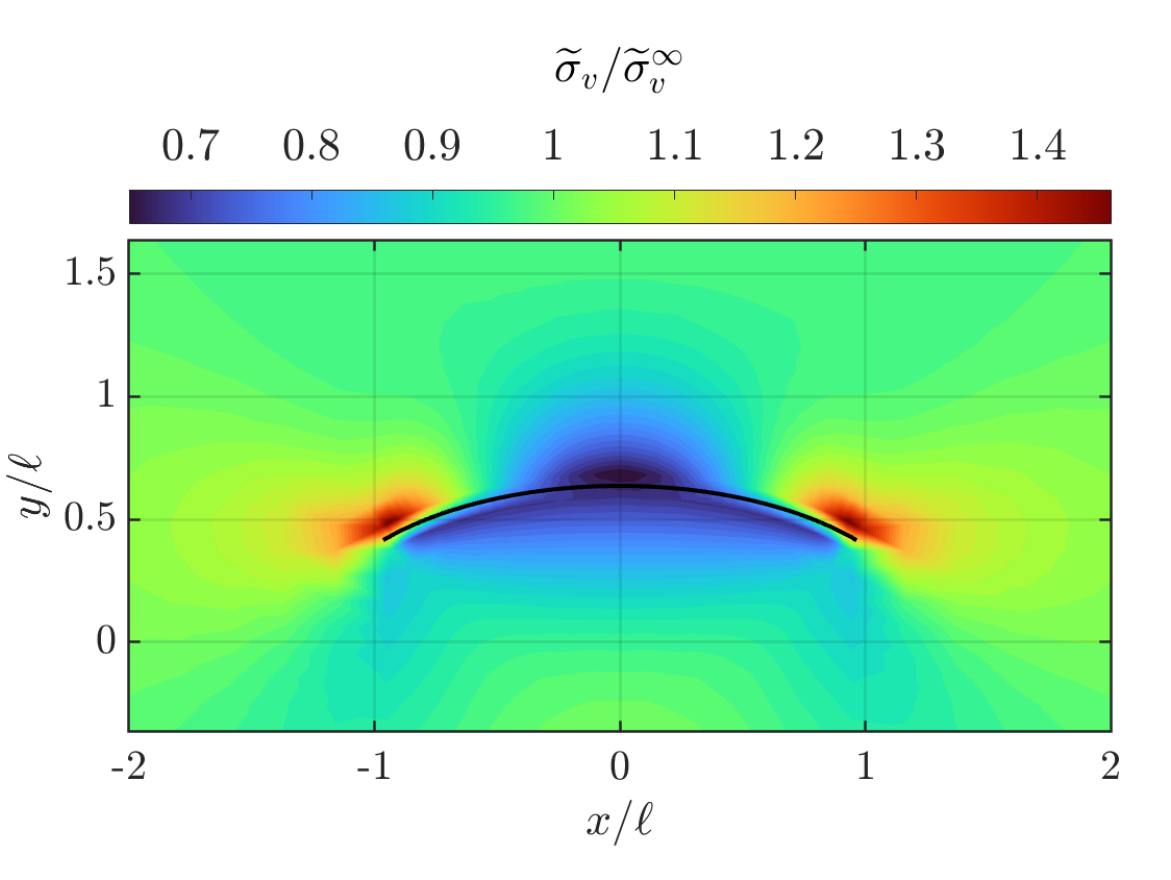}
        \caption{Case (ii)}
        \label{Fig:12.b}
    \end{subfigure}
    \begin{subfigure}{0.495\textwidth}
        \centering
        \includegraphics[width=\linewidth]{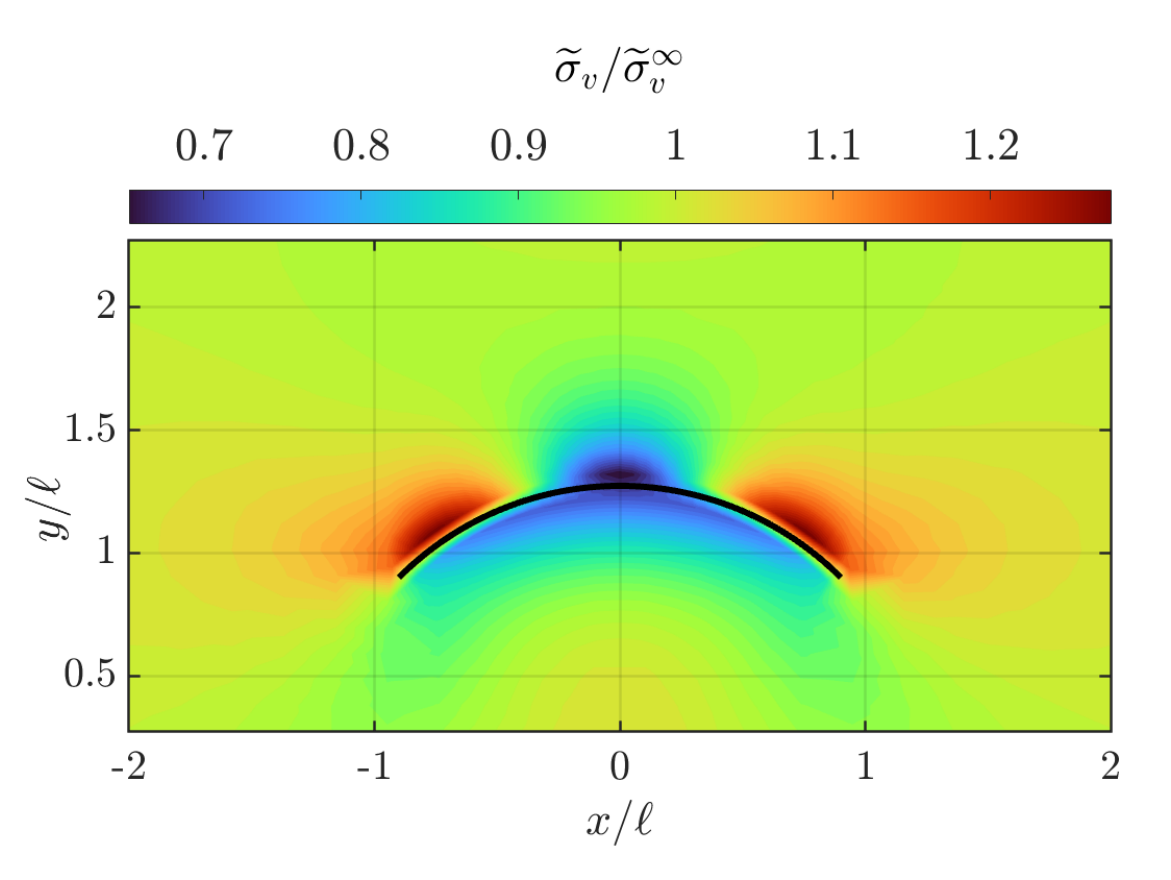}
        \caption{Case (iii)}
        \label{Fig:12.c}
    \end{subfigure}
    \begin{subfigure}{0.495\textwidth}
        \centering
        \includegraphics[width=\linewidth]{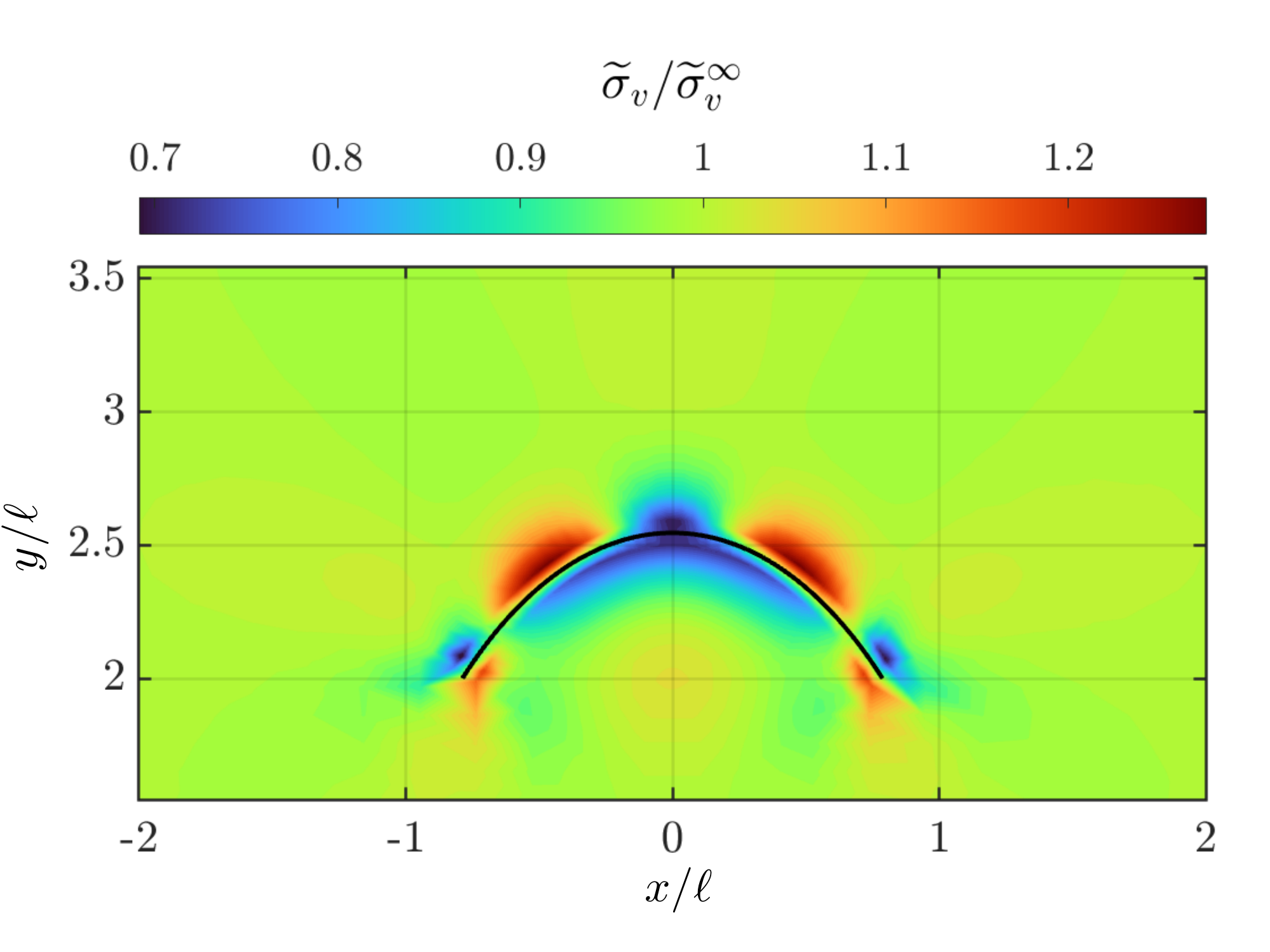}
        \caption{Case (iv)}
        \label{Fig:12.d}
    \end{subfigure}
    \caption{Contours of relative dimensionless Von-Mises stresses for the four surfaces under study.}
    \label{Fig:12}
\end{figure}

The contour plots in Fig.~\ref{Fig:12} highlight the contrasting stress concentration and relaxation patterns associated with different surface geometries. For case (i), the dimensionless Von-Mises stress distribution is symmetric with respect to the $x$- and $y$-axes, while the symmetry with respect to the $y$-axis breaks down for cases (ii)-(iv), reflecting the directional nature of the arc geometry. The contour plots for case (i) feature the wide relaxation zone in the bulk with prominent areas of stress concentrations in relatively narrow vicinities of the tips. The influence of the surface on the bulk material, located right above and below it, rapidly diminishes. In case (ii), the relaxation zone starts to shrink and the areas of stress concentrations start to be more noticeable along the surface boundary. The areas of stress relaxations and concentrations for cases (iii)-(iv) become progressively more localized, indicating that higher curvature variations not only smooth stress distributions, but also confine the perturbations in the bulk material to the areas in the narrower vicinities of the surfaces.

\section{Conclusions}\label{Section:7}

In this paper, we develop a novel, efficient, and robust IGABEM algorithm that could be used to model materials reinforced by ultrathin platelets or sheets.  The Gurtin–Murdoch surface elasticity theory is used in which reinforcement is treated as a material surface located along a sufficiently smooth curve of varying curvature. The integral representations for the elastic fields everywhere in the material system are exact. The unknown components of the surface stress tensor involved in these representations are approximated using the same parametric functions that define the geometry of the curve, namely the NURBS basis functions. This ensures consistency between approximations for the geometry and elastic field and provides the smoothness conditions that are required for the functions involved in the integral representations. 

The proposed algorithm is validated using two benchmark solutions: one involving a surface with an infinite radius of curvature (along a straight segment) and the other involving a surface of a finite but constant radius of curvature (along a circular arc). In both cases, the results obtained with the IGABEM algorithm demonstrated excellent agreement with the benchmark solutions. However, in the latter solutions, the unknowns were approximated globally using various series expansions multiplied by square-root weight functions to enforce the tip conditions. With such approximations, spectral filtering techniques were required for some values of the dimensionless governing parameters to properly capture the near-tip behavior of the fields. In contrast, the IGABEM algorithm does not require filtering and allows for accurate capturing the near-tip effects by using appropriate meshing.

We also report novel results for the problems involving surfaces located along the elliptical curves and, in order to showcase the effects of curvature variations, compare them with those for the problems with surfaces of constant and infinite curvature. 
We demonstrate that the proposed IGABEM algorithm exhibits convergence for both components of the surface stress tensor. 
The contour plots of the relative dimensionless Von Mises stresses presented here reveal a strong influence of curvature variations on the local fields. All these novel results could be used as benchmark solutions for future investigators.

We plan on extending this work in several directions. First, multiple surfaces along arbitrary, sufficiently smooth open and closed curves will be considered, in order to investigate their mutual interactions. Second, homogenization approaches will be used to quantify the influence of curvature and surface tension on the effective mechanical properties of composites. Third, problems with even stiffer reinforcements will be studied using the Steigmann–Ogden surface elasticity theory, which allows to account for bending effects. Finally, the formulation can be generalized to three dimensions, as similar exact integral representations are available in the 3D setting.

\section*{Acknowledgements}\label{Sec:Acknowledgement}
The first and third authors (R.S. Patil and S.G. Mogilevskaya) gratefully acknowledge the support from the U.S. National Science Foundation, award number NSF CMMI -2112894.  R. S. Patil also acknowledges the support from the Hsiao Shaw-Lundquist Fellowship, University of Minnesota.
The second author (Z. Han) gratefully acknowledges the support provided by the Natural Science Foundation of Shanghai Basic Research Funding 25ZR1402003 and by the Fundamental Research Funds for the Central Universities 25D111907.

\appendix
\section{Boundary integral representations}\label{Section:Appendix A}
The displacements components in the Cartesian coordinates can be derived in terms of the local traction jumps by employing Eqs.\;(\ref{Eq:6})-(\ref{Eq:7}) along with the standard coordinate transformation for traction jumps, given by $\Delta t_{1} = \Delta t_{l}\cos{\alpha} + \Delta t_{n}\sin{\alpha}$ and $\Delta t_{2} = \Delta t_{l}\sin{\alpha} - \Delta t_{n}\cos{\alpha}$, where $\alpha = \alpha(s)$ is the angle between the axis $Ox_{1}$ and the tangent at the point $\mathbf{y} \in L$.

That leads to the following expressions:
\begin{align}\label{Eq: Appendix A.1}
    u_{1}(\mathbf{x}) =\ & u_{1}^{\infty}(\mathbf{x}) + \frac{-\kappa}{2\pi\mu(1+\kappa)} \int_{L} \left( \cos{\alpha}\:\Delta t_{l} + \sin{\alpha}\:\Delta t_{n} \right) \ln{r}\, \mathrm{d}s_{\mathbf{y}} \notag \\
    & + \frac{1}{2\pi\mu(1+\kappa)} \int_{L} \left( \cos{\alpha}\:\Delta t_{l} + \sin{\alpha}\:\Delta t_{n} \right) \frac{r_{1}^2}{r^2}\, \mathrm{d}s_{\mathbf{y}} \notag \\ 
    & + \frac{1}{2\pi\mu(1+\kappa)} \int_{L} \left( \cos{\alpha}\:\Delta t_{l} + \sin{\alpha}\:\Delta t_{n} \right) \frac{r_{1}r_{2}}{r^2}\, \mathrm{d}s_{\mathbf{y}} \:,
\end{align}

\begin{align}\label{Eq: Appendix A.2}
    u_{2}(\mathbf{x}) =\ & u_{2}^{\infty}(\mathbf{x})+ \frac{-\kappa}{2\pi\mu(1+\kappa)} \int_{L} \left( \sin{\alpha}\:\Delta t_{l} - \cos{\alpha}\:\Delta t_{n} \right) \ln{r}\, \mathrm{d}s_{\mathbf{y}} \notag \\
    & + \frac{1}{2\pi\mu(1+\kappa)} \int_{L} \left( \sin{\alpha}\:\Delta t_{l} - \cos{\alpha}\:\Delta t_{n} \right) \frac{r_{2}^2}{r^2}\, \mathrm{d}s_{\mathbf{y}} \notag \\ 
    & + \frac{1}{2\pi\mu(1+\kappa)} \int_{L} \left( \cos{\alpha}\:\Delta t_{l} + \sin{\alpha}\:\Delta t_{n} \right) \frac{r_{1}r_{2}}{r^2}\, \mathrm{d}s_{\mathbf{y}} \:,
\end{align}
where, $u_{1}^{\infty}$ and $u_{2}^{\infty}$ are 
\begin{equation}\label{Eq: Appendix A.3}
    \begin{gathered}
        u_{1}^{\infty}(\mathbf{x}) = \left[\frac{(\kappa + 1)\sigma_{11}^{\infty} + (\kappa - 3)\sigma_{22}^{\infty}}{8\mu}\right]x_{1} + \frac{\sigma_{12}^{\infty}}{2\mu}x_{2}\:, \\[10pt]
        u_{2}^{\infty}(\mathbf{x}) =  \frac{\sigma_{12}^{\infty}}{2\mu}x_{1} + \left[\frac{(\kappa - 3)\sigma_{11}^{\infty} + (\kappa + 1)\sigma_{22}^{\infty}}{8\mu}\right]x_{2}\:.
    \end{gathered}
\end{equation}

Similarly, using the coordinate transformation again, the displacement components in the local coordinate system are obtained as $u_{l} = u_{1}\cos{\alpha} + u_{2}\sin{\alpha}$, $u_{n} = u_{1}\sin{\alpha} - u_{2}\cos{\alpha}$. Furthermore, by employing the chain rule for differentiation, we obtain that
\begin{equation}\label{Eq: Appendix A.4}
    \begin{gathered}
     -\frac{u_{l}}{R} + \frac{\partial u_{n}}{\partial s} = \sin{\alpha} \frac{\partial u_{1}}{\partial s} - \cos{\alpha} \frac{\partial u_{2}}{\partial s}\:, \\[10pt]
     \frac{u_{n}}{R} + \frac{\partial u_{l}}{\partial s} =  \sin{\alpha} \frac{\partial u_{2}}{\partial s} + \cos{\alpha} \frac{\partial u_{1}}{\partial s}\:.
    \end{gathered}
\end{equation}

As displacements are continuous when $\mathbf{x} \to \mathbf{y}_0\in L$, one can  obtain the BIEs for $\sigma^{S}(\mathbf{y}_{0})$ and $\omega^{S}(\mathbf{y}_{0})$ by combining Eqs.\;(\ref{Eq:3})-(\ref{Eq:4}) with Eqs.\;(\ref{Eq: Appendix A.1}) to (\ref{Eq: Appendix A.4}) for $\mathbf{x}= \mathbf{y}_0$. The resulting system of singular integral equations in terms of $\sigma^{S},\, \omega^{S}$ and their derivatives has the following forms:
\begin{align}
\sigma^{S}(\mathbf{y}_{0}) &= 
\sigma_{0}(\mathbf{y}_{0})
+ (2\mu_{S} + \lambda_{S})
\Big[
    - \sin \beta(\mathbf{y}_{0}) 
     \Sigma_{1}(\mathbf{y}_{0}) 
    + 
    \cos \beta(\mathbf{y}_{0}) 
    \Sigma_{2}(\mathbf{y}_{0})
\Big] \notag \\[6pt]
&\quad
+ \frac{\kappa (2\mu_{S} + \lambda_{S})}{2\pi \mu (1+\kappa)}
\, \mathrm{C.P.V.} \int_{L}
\Big[
    \sin \beta(\mathbf{y}_{0})\, \textsl{g}_{1}(\mathbf{y})
    - 
    \cos \beta(\mathbf{y}_{0})\, \textsl{g}_{2}(\mathbf{y})
\Big]
\phi_{1}(\mathbf{y}_{0}, \mathbf{y})\, \mathrm{d}s_\mathbf{y}
\notag \\[6pt]
&\quad
- \frac{(2\mu_{S} + \lambda_{S})}{2\pi \mu (1+\kappa)}
\, \mathrm{C.P.V.} \int_{L}
\Big[
    \sin \beta(\mathbf{y}_{0})\, \textsl{g}_{1}(\mathbf{y})
    + 
    \cos \beta(\mathbf{y}_{0})\, \textsl{g}_{2}(\mathbf{y})
\Big]
\phi_{2}(\mathbf{y}_{0}, \mathbf{y})\, \mathrm{d}s_\mathbf{y}
\notag \\[6pt]
&\quad
+ \frac{(2\mu_{S} + \lambda_{S})}{2\pi \mu (1+\kappa)}
\, \mathrm{C.P.V.} \int_{L}
\Big[
    \cos \beta(\mathbf{y}_{0})\, \textsl{g}_{1}(\mathbf{y})
    - 
    \sin \beta(\mathbf{y}_{0})\, \textsl{g}_{2}(\mathbf{y})
\Big]
\phi_{3}(\mathbf{y}_{0}, \mathbf{y})\, \mathrm{d}s_\mathbf{y}
\end{align}

\begin{align}
\omega^{S}(\mathbf{y}_{0})
&=
\cos \beta(\mathbf{y}_{0}) 
\Sigma_{1}(\mathbf{y}_{0})
+
\sin \beta(\mathbf{y}_{0}) 
\Sigma_{2}(\mathbf{y}_{0})
\notag \\[6pt]
&\quad
-\frac{\kappa}{2\pi \mu (1+\kappa)}
\, \mathrm{C.P.V.} \int_{L}
\Big[
    \cos \beta(\mathbf{y}_{0})\, \textsl{g}_{1}(\mathbf{y})
    + 
    \sin \beta(\mathbf{y}_{0})\, \textsl{g}_{2}(\mathbf{y})
\Big]
\phi_{1}(\mathbf{y}_{0}, \mathbf{y})\, \mathrm{d}s_\mathbf{y}
\notag \\[6pt]
&\quad
+\frac{1}{2\pi \mu (1+\kappa)}
\, \mathrm{C.P.V.} \int_{L}
\Big[
    \cos \beta(\mathbf{y}_{0})\, \textsl{g}_{1}(\mathbf{y})
    - 
    \sin \beta(\mathbf{y}_{0})\, \textsl{g}_{2}(\mathbf{y})
\Big]
\phi_{2}(\mathbf{y}_{0}, \mathbf{y})\, \mathrm{d}s_\mathbf{y}
\notag \\[6pt]
&\quad
+\frac{1}{2\pi \mu (1+\kappa)}
\, \mathrm{C.P.V.} \int_{L}
\Big[
    \sin \beta(\mathbf{y}_{0})\, \textsl{g}_{1}(\mathbf{y})
    + 
    \cos \beta(\mathbf{y}_{0})\, \textsl{g}_{2}(\mathbf{y})
\Big]
\phi_{3}(\mathbf{y}_{0}, \mathbf{y})\, \mathrm{d}s_\mathbf{y}
\end{align}
where $\beta(\mathbf{y})$($\beta(\mathbf{y_0}))$ is the angle between the axis $Ox_{1}$ and the normal to the surface $L$ at the point $\mathbf{y}$ ($\mathbf{y_0}$),  $\beta(\mathbf{y})= \alpha(\mathbf{y})- \pi/2$, and 
\begin{align}
\textsl{g}_{1}(\mathbf{y}) 
&= 
-\sin\beta(\mathbf{y})\, \left(\frac{\partial \sigma^{S}}{\partial s} + \sigma_{0} \frac{\omega^{S}}{R}\right)
+ \cos\beta(\mathbf{y})\,\left( - \frac{\sigma^{S}}{R} + \sigma_{0} \frac{\partial\omega^{S}}{\partial s}\right)\,, 
\\[6pt]
\textsl{g}_{2}(\mathbf{y}) 
&= 
\cos\beta(\mathbf{y})\, \left(\frac{\partial \sigma^{S}}{\partial s} + \sigma_{0} \frac{\omega^{S}}{R}\right)
+ \sin\beta(\mathbf{y})\, \left( - \frac{\sigma^{S}}{R} + \sigma_{0} \frac{\partial\omega^{S}}{\partial s}\right)\,,
\end{align}

\begin{align}
\Sigma_{1}(\mathbf{y}_{0})
&=
-\Bigg[
    \frac{(\kappa + 1)\, \sigma_{11}^{\infty} + (\kappa - 3)\, \sigma_{22}^{\infty}}{8\mu}
\Bigg]
\sin \beta(\mathbf{y}_{0})
+
\Bigg[
    \frac{\sigma_{12}^{\infty}}{2\mu}
\Bigg]
\cos \beta(\mathbf{y}_{0})\,,
\\[10pt]
\Sigma_{2}(\mathbf{y}_{0})
&=
-\Bigg[
    \frac{\sigma_{12}^{\infty}}{2\mu}
\Bigg]
\sin \beta(\mathbf{y}_{0})
+
\Bigg[
    \frac{(\kappa - 3)\, \sigma_{11}^{\infty} + (\kappa + 1)\, \sigma_{22}^{\infty}}{8\mu}
\Bigg]
\cos \beta(\mathbf{y}_{0})\,,
\end{align}

\begin{align}
\phi_{1}(\mathbf{y}_{0}, \mathbf{y})
&=
\frac{
    -r_{1}(\mathbf{y}_{0}, \mathbf{y}) \sin \beta(\mathbf{y}_{0})
    + r_{2}(\mathbf{y}_{0}, \mathbf{y}) \cos \beta(\mathbf{y}_{0})
}{
     r^{2}(\mathbf{y}_{0}, \mathbf{y})
},
\\[8pt]
\phi_{2}(\mathbf{y}_{0}, \mathbf{y})
&=
\frac{
    2\, r_{1}(\mathbf{y}_{0}, \mathbf{y})\, r_{2}(\mathbf{y}_{0}, \mathbf{y})
}{
     r^{4}(\mathbf{y}_{0}, \mathbf{y}) 
}
\Big[
    -r_{2}(\mathbf{y}_{0}, \mathbf{y}) \sin \beta(\mathbf{y}_{0})
    - r_{1}(\mathbf{y}_{0}, \mathbf{y}) \cos \beta(\mathbf{y}_{0})
\Big],
\\[8pt]
\phi_{3}(\mathbf{y}_{0}, \mathbf{y})
&=
-\Bigg[
    \frac{
        r_{2}(\mathbf{y}_{0}, \mathbf{y})
    }{
        r^{2}(\mathbf{y}_{0}, \mathbf{y})
    }
    - 
    \frac{
        2\, r_{1}^{2}(\mathbf{y}_{0}, \mathbf{y})\, r_{2}(\mathbf{y}_{0}, \mathbf{y})
    }{
        r^{4}(\mathbf{y}_{0}, \mathbf{y})
    }
\Bigg]
\sin \beta(\mathbf{y}_{0})
\nonumber \\[4pt]
&\quad
+\Bigg[
    \frac{
        r_{1}(\mathbf{y}_{0}, \mathbf{y})
    }{
        r^{2}(\mathbf{y}_{0}, \mathbf{y})
    }
    - 
    \frac{
        2\, r_{1}(\mathbf{y}_{0}, \mathbf{y})\, r_{2}^{2}(\mathbf{y}_{0}, \mathbf{y})
    }{
        r^{4}(\mathbf{y}_{0}, \mathbf{y})
    }
\Bigg]
\cos \beta(\mathbf{y}_{0}).
\end{align}

\section{Convergence study}\label{Section:Appendix B}
For the convergence study, two surface geometries are considered: (i) a straight segment and (ii) a circular arc. The material properties of the matrix and surfaces, as well as the loading conditions, are taken as described in Section \ref{Sub-Section:6.1}. The relative $\mathrm{L}^{2}$-norm errors for surface stress components are defined as 
\begin{equation}
E_{\widetilde{\sigma}^{S}} =
\sqrt{ \frac{\displaystyle\sum_{k=1}^{N_{\mathrm{p}}}
\left[ e_{\widetilde{\sigma}^{S}}(\xi_{k}) \right]^2}{ \displaystyle\sum_{k=1}^{N_{\mathrm{p}}}
\left[ \widetilde{\sigma}^{S}_{\mathrm{ref}}(\xi_{k}) \right]^2}}\,,\qquad
E_{\omega^{S}} =
\sqrt{
\frac{
\displaystyle\sum_{k=1}^{N_{\mathrm{p}}}
\left[ e_{\omega^{S}}(\xi_{k}) \right]^2
}{
\displaystyle\sum_{k=1}^{N_{\mathrm{p}}}
\left[ \omega^{S}_{\mathrm{ref}}(\xi_{k}) \right]^2
}
}\,,
\end{equation}
where $N_{\mathrm{p}} = 200$ denotes the total number of uniformly distributed sampling points $\xi_{k}$ in the parametric space and
\begin{equation}
  e_{\widetilde{\sigma}^{S}}(\xi_{k}) = \widetilde{\sigma}^{S}(\xi_{k}) - \widetilde{\sigma}^{S}_{\mathrm{ref}}(\xi_{k})\,,\qquad 
  e_{\omega^{S}}(\xi_k) = \omega^{S}(\xi_k) - \omega^{S}_{\mathrm{ref}}(\xi_k)\,.
\end{equation}
\noindent

Here, $\widetilde{\sigma}^{S}(\xi_k)$ and $\widetilde{\sigma}^{S}_{\mathrm{ref}}(\xi_k)$ represent the normalized first surface stresses component obtained from the coarser and reference boundary meshes, respectively, while $\omega^{S}(\xi_k)$ and $\omega^{S}_{\mathrm{ref}}(\xi_k)$ denote those for the second surface stress component. Since analytical solutions are not available for the two cases, the results obtained using 102 degrees of freedom (control points) are taken as reference solutions $\widetilde{\sigma}^{S}_{\mathrm{ref}}$ and $\omega^{S}_{\mathrm{ref}}$. Thus, $e_{\widetilde{\sigma}^{S}}(\xi_k)$ and $e_{\omega^{S}}(\xi_k)$ quantify the differences between the coarser and reference solutions at the point $\xi_k$.

The relative $\mathrm{L}^{2}$-norm errors as functions of the degrees of freedom (i.e., the number of control points given by $N_{e}+2$) are presented in Fig.~\ref{Fig:Appendix B.1}, with results for the surface along straight segment shown in Fig.~\ref{Fig:Appendix B.1a} and for along circular arc in Fig.~\ref{Fig:Appendix B.1b}.
\setcounter{figure}{0}
\renewcommand{\thefigure}{B.\arabic{figure}}
\begin{figure}[H]
  \centering
  \begin{subfigure}[b]{0.49\textwidth}
      \centering
      \includegraphics[angle=0, width=1\textwidth]{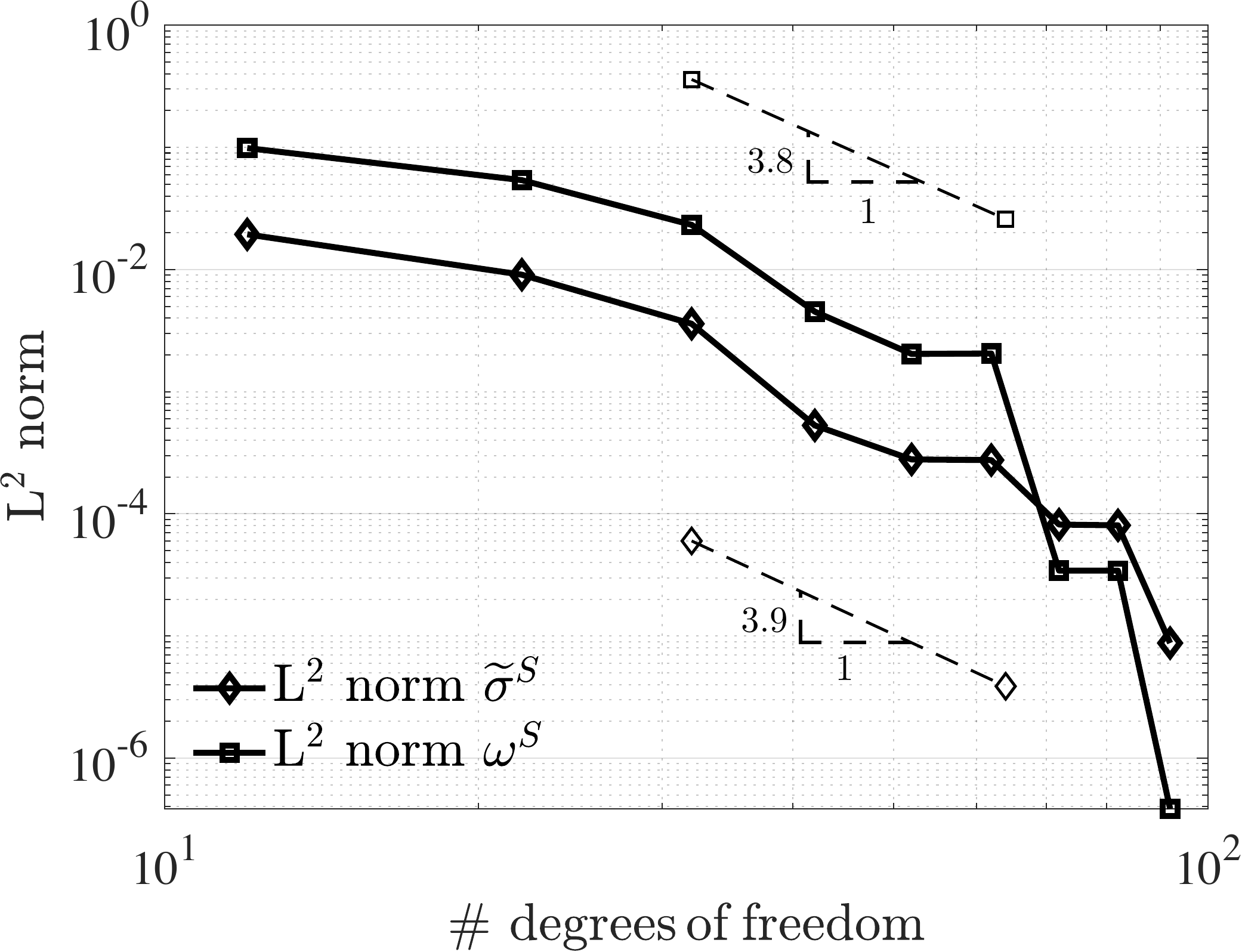}
      \subcaption{}
      \label{Fig:Appendix B.1a}
  \end{subfigure}
  \hfill
  \begin{subfigure}[b]{0.49\textwidth}
      \centering
      \includegraphics[angle=0, width=1\textwidth]{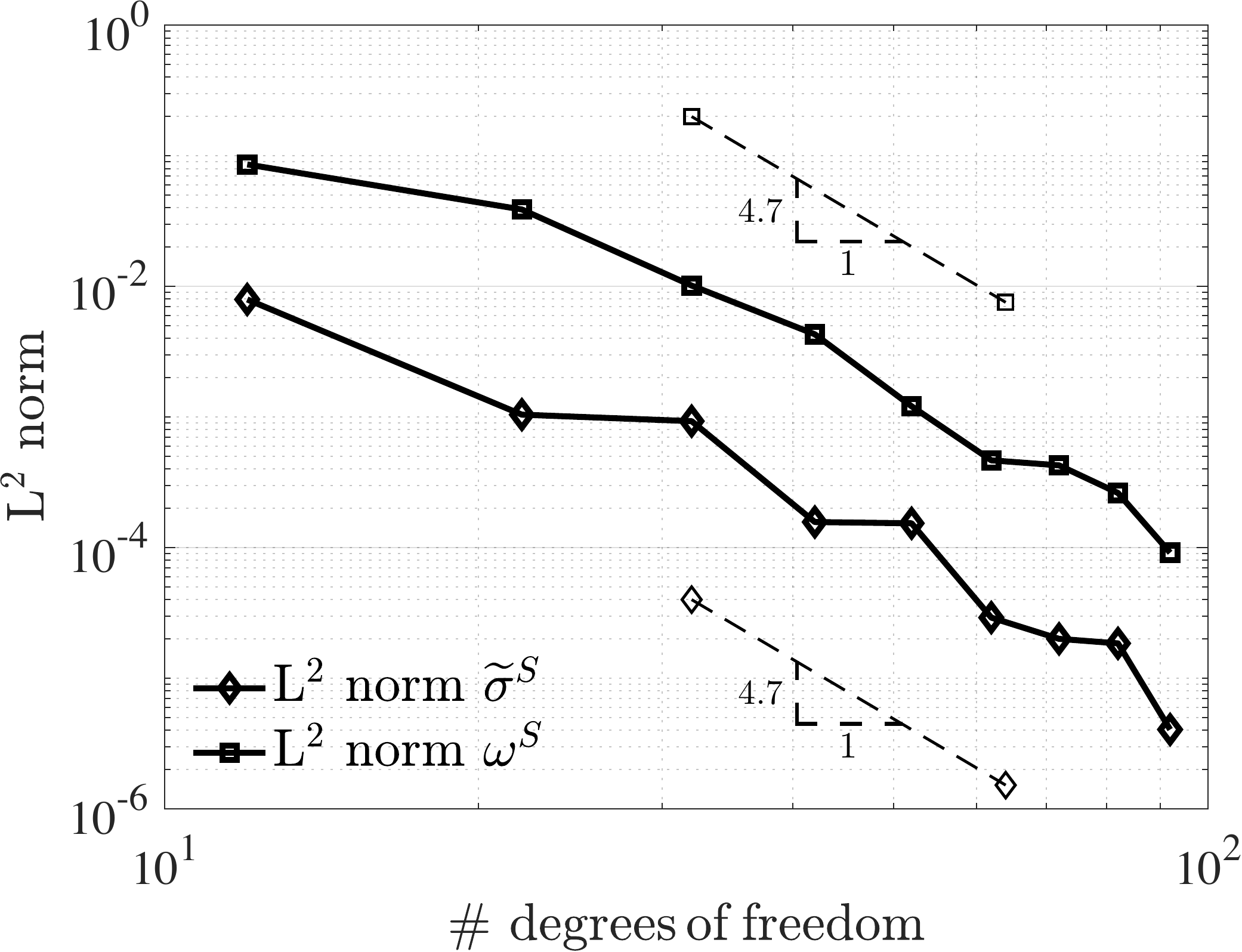}
      \subcaption{}
      \label{Fig:Appendix B.1b}
  \end{subfigure}
  \caption{Convergence in the $\mathrm{L}^{2}$ norm error with quadratic elements for surface along a (a) straight line (left) and (b) circular arc (right).}
  \label{Fig:Appendix B.1}
\end{figure}

\begin{figure}[H]
  \centering
  \begin{subfigure}[b]{0.49\textwidth}
      \centering
      \includegraphics[angle=0, width=1\textwidth]{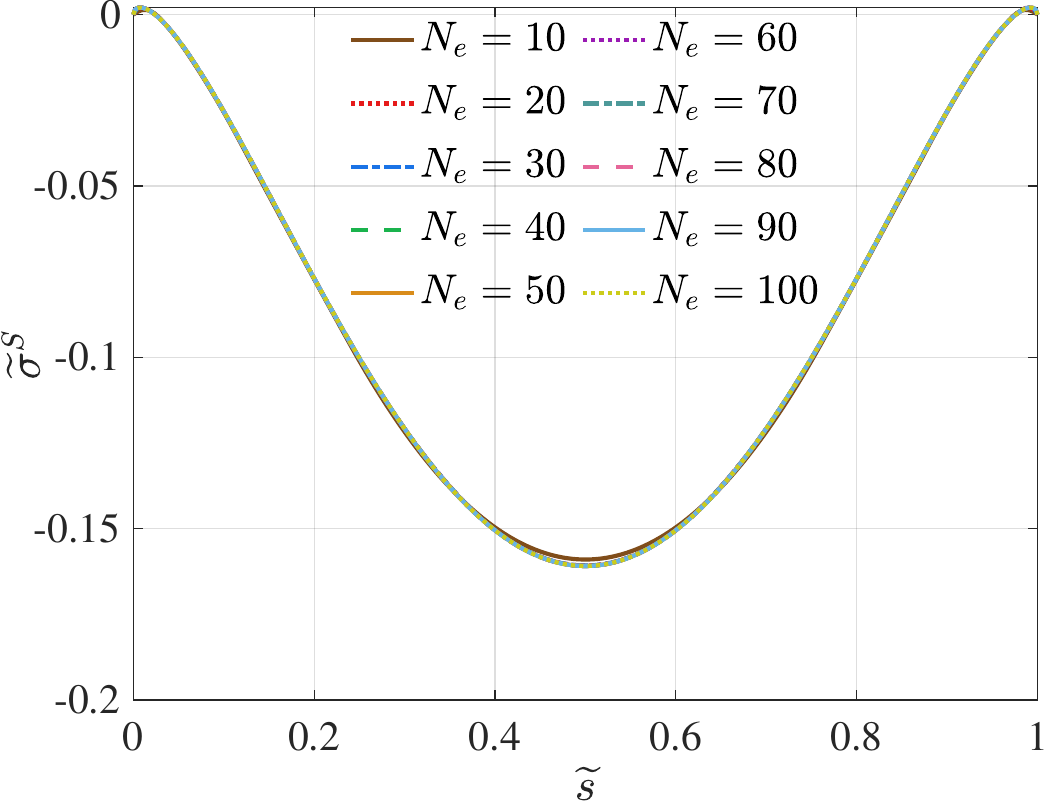}
      \subcaption{}
      \label{Fig:Appendix B.2a}
  \end{subfigure}
  \hfill
  \begin{subfigure}[b]{0.49\textwidth}
      \centering
      \includegraphics[angle=0, width=1\textwidth]{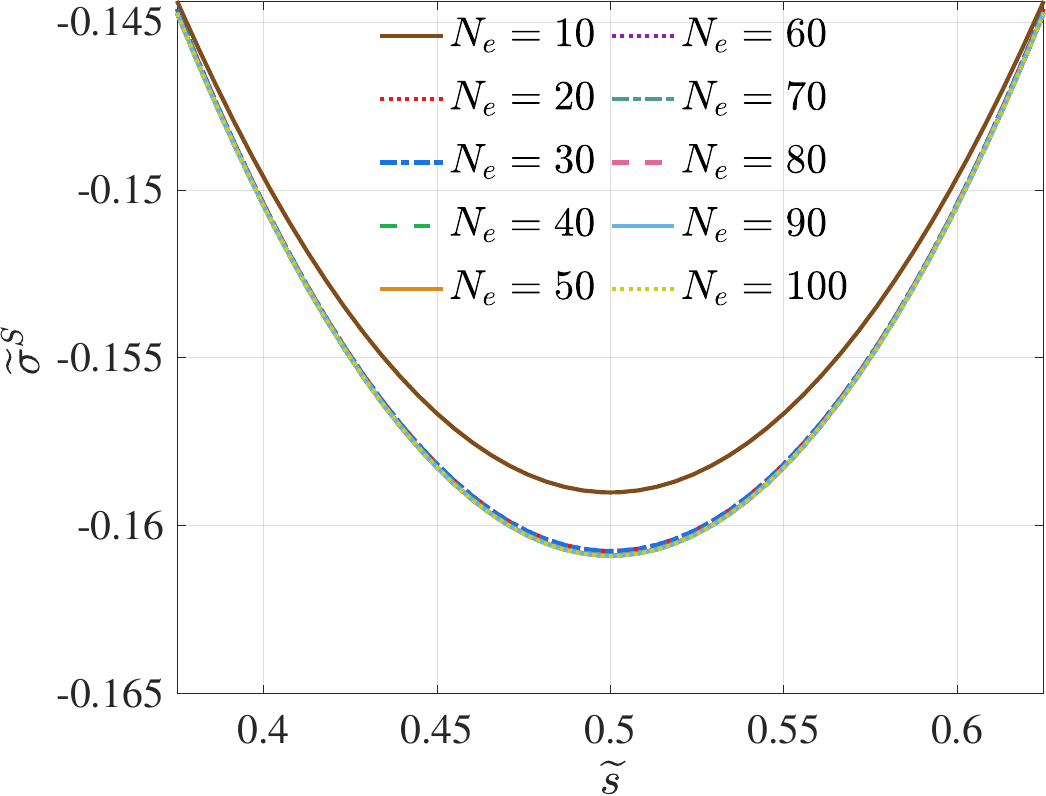}
      \subcaption{}
      \label{Fig:Appendix B.2b}
  \end{subfigure}
  \caption{(a) The plots of $\widetilde{\sigma}^{S}$ for surface along circular arc for different values of $N_{e}$ and (b) zoomed-in view highlighting the detailed convergence behavior.}
  \label{Fig:Appendix B.2}
\end{figure}

\begin{figure}[H]
  \centering
  \begin{subfigure}[b]{0.49\textwidth}
      \centering
      \includegraphics[angle=0, width=1\textwidth]{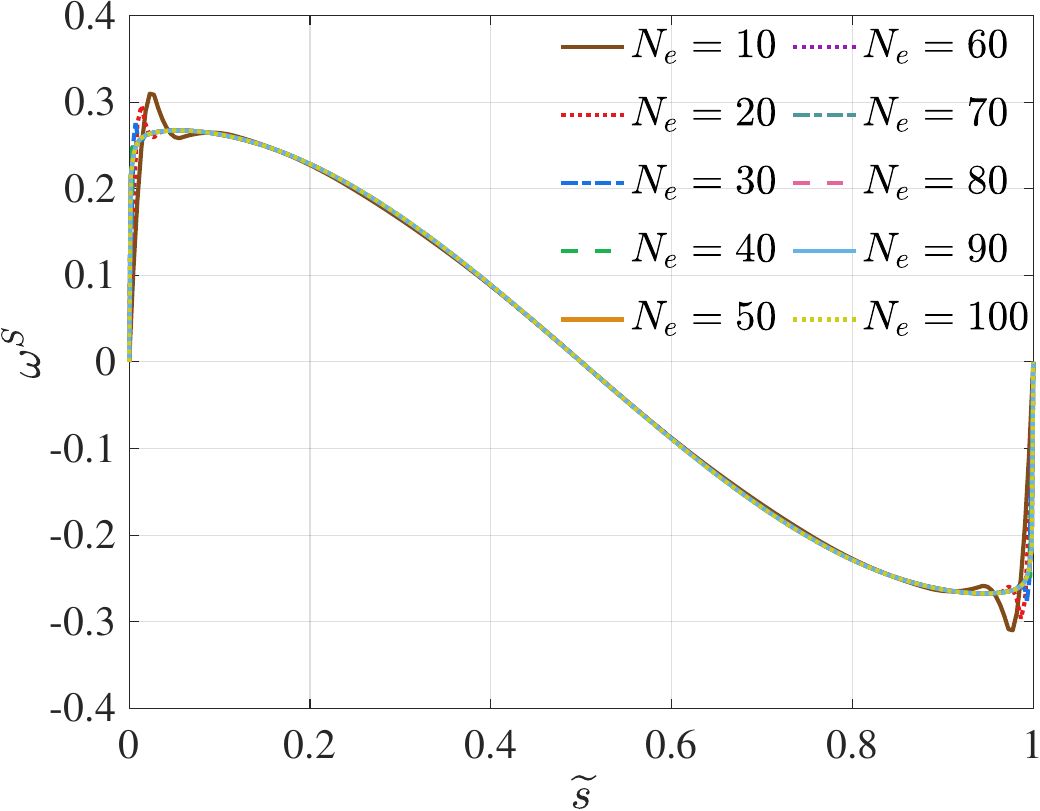}
      \subcaption{}
      \label{Fig:Appendix B.3a}
  \end{subfigure}
  \hfill
  \begin{subfigure}[b]{0.49\textwidth}
      \centering
      \includegraphics[angle=0, width=1\textwidth]{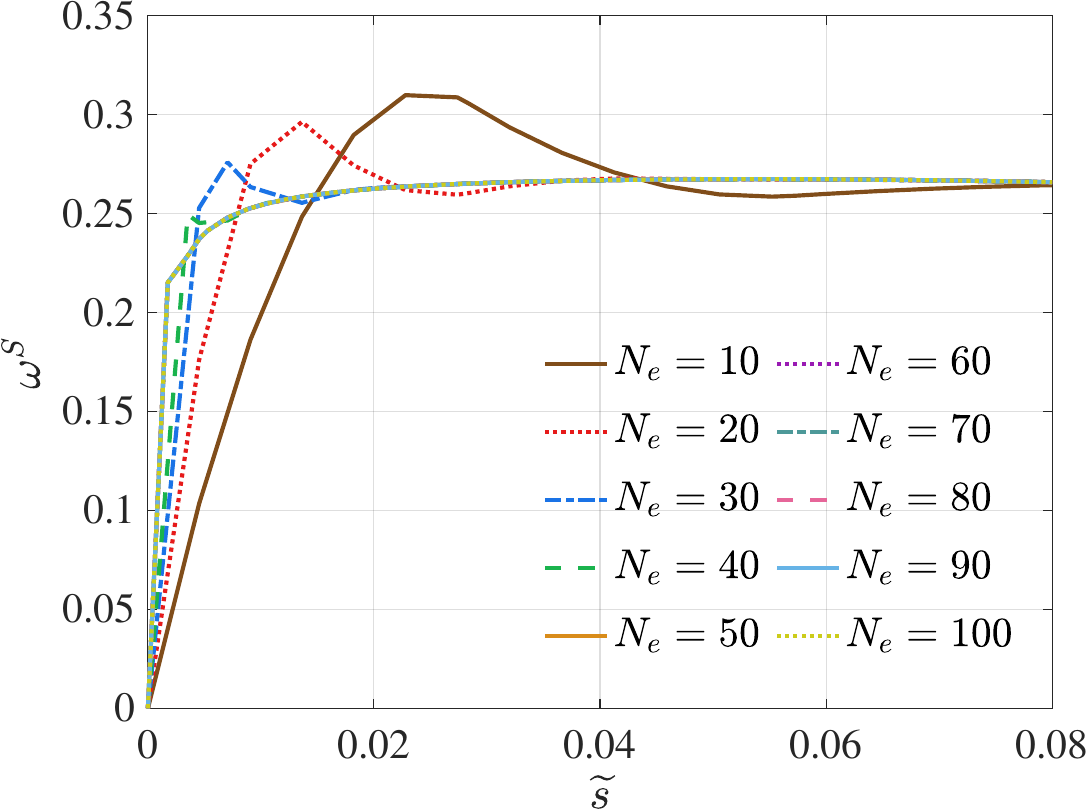}
      \subcaption{}
      \label{Fig:Appendix B.3b}
  \end{subfigure}
  \caption{(a) The plots of $\omega^{S}$ for surface along circular arc for different values of $N_{e}$ and (b) zoomed-in view highlighting the detailed convergence behavior.}
  \label{Fig:Appendix B.3}
\end{figure}

As can be seen in Fig.~\ref{Fig:Appendix B.1}, the $\mathrm{L}^{2}$ norm errors rapidly decrease with increasing number of degrees of freedom. The error for $\omega^{S}$ is consistently higher than that for $\widetilde{\sigma}^{S}$, reflecting a higher sensitivity of the components to discretizations. In addition, to further illustrate convergence, we plotted in Figs.~\ref{Fig:Appendix B.2} and \ref{Fig:Appendix B.3} the distributions of $\widetilde{\sigma}^{S}$ and $\omega^{S}$  for the surface along the circular arc as functions of $N_{e}$, ranging from 10 to 100.

From Fig.~\ref{Fig:Appendix B.2} it can be concluded that the results for $\widetilde{\sigma}^{S}$ converge fast with refinement, and the plot for $N_{e} \geq 40$ is practically indistinguishable from that for $N_{e}=100$. The zoomed-in view further highlights the convergence as the mesh is refined. Similar conclusions for all $N_{e}$ can be drawn from the analysis of the plots of $\omega^{S}$ in Fig.~\ref{Fig:Appendix B.3}.  Most discrepancies in the results obtained with different meshes occur at the small intervals near the tips, $\widetilde{s}=0$ and $\widetilde{s}= 1$, as shown in Fig.~\ref{Fig:Appendix B.3b} for the tip $\mathbf{a}$. The plots for coarser meshes ($N_{e}=10–40$) exhibit localized oscillations in these regions, while those oscillations are removed with further refinement ($N_{e} \geq 50$), and the results converge smoothly to zero at the tips. So, it can be concluded that the mesh with $N_{e}=50 $ is sufficient to capture the behavior at the tips with high fidelity.

\bibliographystyle{elsarticle-num}  
 
\bibliography{references_Mendeley_trial}    

\end{document}